\documentclass{amsart}

\usepackage[english]{babel}
\usepackage[utf8]{inputenc}
\usepackage{times}
\usepackage[T1]{fontenc}
\usepackage{amsmath,amsthm,amssymb,xspace}
\usepackage{xcolor}
\usepackage{tikz}
\usepackage{tabularx}
\usepackage[normalem]{ulem}
\usetikzlibrary{decorations.pathmorphing,positioning,patterns}

\ifx\UseOption\undefined
\def\UseOption{pdftikz,dessinempost,color}
\fi
\usepackage{optional}

\opt{color}{
\colorlet{couleurlien}{green!20!black}
\colorlet{couleurlink}{red!60!black}
}
\opt{bw}{
\colorlet{couleurlien}{black}
\colorlet{couleurlink}{black}
}

\usepackage[pdftex,colorlinks,citecolor=couleurlien,linkcolor=couleurlink]{hyperref}
\newcommand\sshift[1]{\ensuremath{\mathbf{#1}}}
\newcommand\setofconf[1]{\ensuremath{\mathbf{#1}}}

\newcommand\plan[1][\groupe{\Z^2}]{\ensuremath{\mathbf{#1}}}

\newcommand\couleur[1]{\mathbf{#1}}
\newcommand\motif[1]{\ensuremath{\mathbf{#1}}}

\newcommand\ensdef[1]{\ensuremath{\mathbf{#1}}}
\newcommand\domain[1]{\ensuremath{\mathcal{D}_{#1}}}

\newcommand\enspav{\mathbf{X}}
\newcommand\dist[1]{\ensuremath{d_{#1}}}
\newcommand\distance[3]{\ensuremath{\dist{#1}(#2,#3)}}
\newcommand\norme[2]{\ensuremath{||#2||_{#1}}}
\newcommand\ouvert[1]{\ensuremath{\mathcal{#1}}}
\newcommand\open[1]{\ensuremath{\ouvert{O_{\motif{#1}}}}}
\newcommand\espace[1][\plan]{\ensuremath{\alphab^{#1}}}

\newcommand\conf[1]{#1}

\newcommand\orbit[1]{\ensuremath{\mathcal{O}\left(\conf{#1}\right)}}

\newcommand\shiftalone[1]{\ensuremath{\sigma_{#1}}}
\newcommand\shift[2]{\ensuremath{\shiftalone{#1}(\conf{#2})}}
\newcommand\class[1]{\ensuremath{\left<#1\right>}}
\newcommand\ForbPat{\ensuremath{\mathcal{F}}}
\newcommand\Tilings[1][\ForbPat]{\ensuremath{\enspav_{#1}}}

\newcommand\rk[1]{\ensuremath{\rho_{\sshift{#1}}}}
\newcommand\engendre[1]{\ensuremath{\overline{\orbit{#1}}}}

\newcommand\Lang[1]{\ensuremath{\mathcal{L}(#1)}}
\newcommand\vdistance[1]{#1}
\newcommand\boule[3]{\ensuremath{\mathcal{B}_{#1}(#2,\vdistance{#3})}}

\newcommand\cardcont{\ensuremath{2^{\aleph_0}}}
\newcommand\epsi{\varepsilon}
\newcommand\arbre[1]{#1}

\newcommand\compl[1]{\ensuremath{#1^C}}
\newcommand\complf[3]{\ensuremath{\mathcal{C}_{#1,#2}(\conf{#3}})}

\newcommand\alphab{\ensuremath{\couleur{\Sigma}}}

\newcommand\cellcouleur[1]{\tikz[scale=.2]{\filldraw[fill=#1] (0,0) rectangle +(1,1);}}
\newcommand\latin[1]{\emph{#1}}

\newcommand\ie{\latin{i.e.}, }
\newcommand\eg{\latin{e.g.}, }
\newcommand\Z{\ensuremath{\mathbb{Z}}}
\newcommand\N{\ensuremath{\mathbb{N}}}

\newcommand\bandeinf[1]{\ensuremath{\motif{#1}^{\updownarrow v}}}

\newcommand\derivee[2]{\ensuremath{\sshift{#1}^{(#2)}}}

\opt{color}{
\newcommand\mgreen{green!75!black!50}

\newcommand\orange{orange}
\newcommand\yellow{yellow}
\newcommand\myellow{yellow!60}
\newcommand\purple{purple}

}
\opt{bw}{
\newcommand\mgreen{black!35}

\newcommand\orange{black!50}
\newcommand\yellow{black!15}
\newcommand\myellow{white}
\newcommand\purple{black!65}

}

\newtheorem{theorem}{Theorem}[section]
\newtheorem{cor}{Corollary}

\newtheorem{lemma}[theorem]{Lemma}
\newtheorem{prop}[theorem]{Proposition}

\theoremstyle{definition}
\newtheorem{defi}[theorem]{Definition}
\newtheorem{op}{Problem}

\theoremstyle{remark}

\newtheorem{claim}{Claim}[section]

\title{Structuring multi-dimensional subshifts}

\begin{document}

\author{Alexis Ballier}
\address{Centro de Modelamiento Matemático,
Universidad de Chile,
Av. Blanco Encalada 2120, Piso 7,
Santiago, Chile}
\email{aballier@dim.uchile.cl}
\thanks{Alexis Ballier is supported by the FONDECYT Postdoctorado Proyecto
3110088}
\author{Emmanuel Jeandel}
\address{LORIA
        Campus Scientifique - BP 239
        54506 VANDOEUVRE-LES-NANCY
FRANCE}
\email{emmanuel.jeandel@loria.fr}
\thanks{Emmanuel Jeandel is supported by ANR-09-BLAN-0164 and PEPS ATIC.}

\thanks{Part of this work has been presented at the conference STACS
2008~\cite{stacsstructuralaspect} and was done when Alexis Ballier was a
PhD Student at Université de Provence.
We therefore thank Bruno Durand for his supervision and his work on this
paper~\cite{stacsstructuralaspect} that are the foundations of the
results presented here.}


\begin{abstract}
We study two relations on multi-dimensional subshifts: A pre-order based on the
patterns configurations contain and the Cantor-Bendixson rank.
We exhibit several structural properties of two-dimensional subshifts: We
characterize the simplest aperiodic configurations in countable SFTs, 
we give a combinatorial characterization of uncountable subshifts,
we prove that there always exists configurations without any periodicity but
that have the simplest possible combinatorics in countable SFTs.
Finally, we prove that some Cantor-Bendixson ranks are impossible for countable
SFTs, leaving only a few unknown cases.
\end{abstract}

\maketitle


\section*{Introduction}

The most famous formalism for \emph{tilings of the discrete plane} has been
introduced by Hao Wang in order to study the decidability of a given class of
first order
formulae~\cite{wangpatternrecoII,wangpatternrecoI,wangdominoesaea,BorgerGG1997}.
This formalism is now known as \emph{Wang tiles} and is very simple: Consider a
finite set of square unit tiles with colors on their borders and allow to put
two such tiles side by side only if their adjacent colors match.
Because of the simplicity of the model, Wang conjectured that any set of tiles
that can be used to fill the plane can tile the plane in a periodic
way~\cite{wangpatternrecoII}, which would give a decision procedure for the
\emph{domino problem}: deciding if a set of tiles can tile the plane.
Robert Berger proved a couple of years later that there exist \emph{aperiodic
sets of tiles}~\cite{berger,bergerthesis}: sets of tiles that can fill the plane
but never in a periodic way; he also proved that the domino problem is
undecidable.
A lot of undecidability results for parallel computation models such as
\emph{cellular automata} are reductions to the domino
problem~\cite{revers-2D,revers-2D-comp,hurd92,kari92a,kari-limit-2D,AanderaaL74}.

\emph{Symbolic dynamics} deal with the same objects.
The motivations come from the work of Morse and Hedlund~\cite{SymbDyn} who
wanted to study discretizations of dynamical systems modeled by the \emph{action
of the shift map} over \emph{subshifts}.
In this field, the literature for dimension~$1$ is huge~\cite{marcuslind} and
has recently been extended to higher dimensions~\cite{lind}.
In this paper, we adopt the vocabulary and definitions from symbolic dynamics
since these are the definitions that can be generalized easily: tilings of the
plane by Wang tiles correspond to \emph{subshifts of finite type} in
dimension~$2$ from symbolic dynamics, a tiling is a \emph{configuration} in such
a subshift.
We often do not need the finite type hypothesis and thus can easily state more
general results using this formalism.

Our aim is to study and understand the structural and combinatorial properties
of tilings of the plane.
We want to obtain a precise description of the configurations that can appear in
subshifts (of finite type).
We mainly use two tools:
\begin{itemize}
\item A pre-order, $\preceq$, introduced in \cite{Durand99}, which compares the finite parts of two
configurations
\item The \emph{Cantor-Bendixson} derivation and rank which comes from the
topological properties of subshifts.
\end{itemize}
With these tools we are able to better understand subshifts and sort them in three types:
Independently of the continuum hypothesis, subshifts are either finite,
countably infinite or have the cardinality of continuum.
In particular, infinite countable subshifts of finite type admit a bi-periodic
configuration and a configuration with exactly one direction of periodicity
(Theorem~\ref{thm:sftdnbalorsconf1per} and \cite{stacsstructuralaspect}).
In addition, if such a subshift contains configuration without any periodicity,
Theorem~\ref{thm:lvl2exists} tells us that there exists a simplest (for
$\preceq$) non-periodic configuration.
Theorem~\ref{thm:caracsshiftnondnb} characterizes the cause of uncountability of
subshifts by means of the configurations they contain.

\section{Basic definitions and properties}
\label{sec:defs}

\subsection{Notations}

$\alphab$ is a \emph{finite set} (represented by colors on figures) and called
the \emph{alphabet}.
We focus on dimension~$2$, hence a \emph{configuration} is an element of
$\espace$ or equivalently an application from $\plan$ to $\alphab$.
A \emph{shift} of vector $v \in \Z^2$, denoted by $\shiftalone{v}$, is the
application from
$\espace$ to $\espace$ defined by: $\shift{v}{x}(a)=x(a-v)$.
A configuration $\conf{x}$ is said to be \emph{periodic} if there exists two
independent vectors $v$, $v'$ so that $\shift{v}{x} = \shift{v'}{x} = \conf{x}$.
Equivalently, $\conf{x}$  has a finite orbit under the $\Z^2$ action by shifts.
We say that a set of configurations $\setofconf{X}$ is \emph{shift-invariant} if for
any $v\in\plan$, $\conf{x}$ is in $\setofconf{X}$ whenever $\shift{v}{x}$ is in
$\setofconf{X}$.
A \emph{pattern} $\motif{P}$ is a finite restriction of a configuration, \ie 
an element of $\alphab^V$ for some finite subset $V$ of $\Z^2$.
We denote the domain of definition of $\motif{P}$ by $\domain{P}$.

A pattern \emph{appears} in a configuration $\conf{c}$ (resp.\ in some other pattern 
$\motif{P'}$) if it can be found somewhere in $\conf{c}$ (resp.\ in
$\motif{P'}$); \ie if there exists a vector $v \in \plan$ such that
$\conf{c}(x+v) = \motif{P}(x)$ on the domain of $\motif{P}$ (resp.\ if
$\motif{P'}(x+v)$ is defined for $x\in V$ and $\motif{P'}(x+v) =
\motif{P}(x)$).
In order to simplify the notations, in this article we will
consider patterns up to a shift, that is: Two patterns $\motif{P}$ and
$\motif{P'}$ will be considered equal if $\motif{P}$ appears in $\motif{P'}$ and
$\motif{P'}$ appears in $\motif{P}$.

Given a set of patterns $\ForbPat$, the \emph{subshift} of $\espace$ defined by
forbidding $\ForbPat$, denoted by $\Tilings$, is the set of configurations of
$\espace$ where no pattern of $\ForbPat$ appear.
$\Tilings$ is said to be a \emph{Subshift of Finite Type} (\emph{SFT} in short)
if there exists a finite set of patterns $\ForbPat'$ such that
$\Tilings=\Tilings[\ForbPat']$.
$\Tilings$ is said to be an \emph{effective subshift} if there exists a
recursively enumerable set of patterns $\ForbPat'$ such that
$\Tilings=\Tilings[\ForbPat']$.
A subshift $\sshift{X}$ is said to be \emph{minimal} if it has no proper
subshift in the that that if $\sshift{Y}\subseteq\sshift{X}$ is a subshift then
either $\sshift{Y}=\emptyset$ or $\sshift{Y}=\sshift{X}$.

\subsection{Languages}
As the subshifts are defined based on the patterns they do not contain, it is
natural to consider the following pre-order:

\begin{defi}
\label{def:preorder}
Let $\conf{x}$ and $\conf{y}$ be two configurations, we say that
$\conf{x}\preceq \conf{y}$ if any pattern that appears in $\conf{x}$ also
appears in $\conf{y}$.
\end{defi}

The relation $\preceq$ is a partial pre-order: It is transitive
and symmetric.
We say that two configurations $\conf{x}$ and $\conf{y}$ are equivalent if
$\conf{x}\preceq\conf{y}$ and $\conf{y}\preceq\conf{x}$.  We denote this
relation by $\conf{x}\approx\conf{y}$.  In this case, $\conf{x}$ and $\conf{y}$
contain the same patterns. 
Note that $\conf{x}$ and $\conf{y}$ are equivalent if they are equal up to
shift but the converse is not true in general.
We  write
$\conf{x}\prec\conf{y}$ if $\conf{x}\preceq\conf{y}$ and
$\conf{x}\not\approx\conf{y}$.

Given a configuration $\conf{x}$ we denote by $\Lang{\conf{x}}$, called the
\emph{language of \conf{x}}, the set of all patterns that appear in $\conf{x}$.
By extension, for a set of configurations $\setofconf{X}$, we denote by
$\Lang{\setofconf{X}}$ the set
$\bigcup_{\conf{x}\in\setofconf{X}}\Lang{\conf{x}}$.
By definition $\conf{x}\preceq\conf{y}$ if and only if
$\Lang{\conf{x}}\subseteq\Lang{\conf{y}}$ and a configuration $\conf{x}$ is in a
subshift $\Tilings$ if and only if $\Lang{\conf{x}}\cap\ForbPat=\emptyset$.
Given a configuration $\conf{x}$, the set
$\engendre{x}=\left\{\conf{y}\in\espace, \conf{y}\preceq\conf{x}\right\}$ is a
subshift since it is the subshift $\Tilings[\compl{\Lang{\conf{x}}}]$ where
$\compl{\Lang{\conf{x}}}$ is the set of all patterns that are not in
$\Lang{\conf{x}}$.

\subsection{Topology}

In symbolic dynamics~\cite{gohe,Hedlund69,lind,marcuslind}, subshifts
are seen as topological spaces on which the shift map acts. 
We embed $\alphab$ with the discrete topology and $\espace$ with the product
topology also known as the Cantor topology.
A basis of the topology we obtain is given by the cylinder sets, to every
pattern $\motif{P}$ corresponds a natural open set which is the set of
configurations agreeing with $\motif{P}$ on its domain:
\[
\open{P}=\left\{\conf{c}\in\espace, \conf{c}_{|\domain{P}}=\motif{P}\right\}
\]
Every such $\open{P}$ is also closed and thus $\espace$ is $0-$dimensional as it
has a topology basis of clopen sets.
Since $\alphab$ is finite, it is compact, and by Tychonoff's theorem we obtain
that $\espace$ is also compact.
The same topology can also be obtained by the following metric:
\[
\distance{}{\conf{x}}{\conf{y}}=2^{-\min\left\{\norme{}{v},
\conf{x}(v)\neq\conf{y}(v)\right\}}
\]
Hedlund~\cite{Hedlund69} proved that \emph{subshifts} as defined in the previous
section correspond to \emph{closed and shift-invariant} subsets of $\espace$.
The notation $\engendre{x}$ from the previous section stands for the \emph{orbit
closure} of $\conf{x}$:
$\engendre{x}=\left\{\conf{y},\conf{y}\preceq\conf{x}\right\}=\overline{\left\{\shift{v}{x},
v\in\plan\right\}}$.
SFTs can also be characterized in a more topological way:
\begin{prop}[Folklore, \cite{kschmidtalgzdactions}]
\label{prop:schmidtcaracsft}
A subshift $\sshift{X}$ is of finite type if and only if for any decreasing
chain of subshifts $(\sshift{X}_i)_{i\in\N}$ such that 
$\cap_{i\in\N}\sshift{X}_i=\sshift{X}$ then there exists $N\in\N$ such that
$\sshift{X_N}=\sshift{X}$.
\end{prop}

Proposition~\ref{prop:schmidtcaracsft} can be written combinatorially as follows:
Let $\ForbPat$ be a set of patterns, if $\Tilings$ is an SFT then there exists a
finite subset $\ForbPat'$ of $\ForbPat$ such that $\Tilings=\Tilings[\ForbPat']$.
With this approach, it can easily be proved that when a subshift $\sshift{X}$ is
finite then it is of finite type.

\section{Structural properties}
\label{sec:struct}

We start by some general structural properties of subshifts that are used
thorough this paper.
The first question that arises when trying to understand the pre-order $\preceq$
is the existence of extremal points: It is well known that any subshift contains
a minimal point for $\preceq$~\cite{birkmin,Durand99} and we proved that any
subshift also contains a maximal
point~\cite[Theorem~$3.2$]{stacsstructuralaspect}.
A way to better understand the structure of subshifts with the $\preceq$
pre-order is to try to see this order as a discrete one by generalizing the
notion of minimality to the one of levels:

\begin{defi}[Levels]
\label{defi:levels}
Let $\conf{c}$ be an element of a subshift $\Tilings[]$, we define the level of
$\conf{c}$ by the length of the longest strictly decreasing chain for $\preceq$
in $\Tilings[]$ starting at $\conf{c}$.

If such chains can be arbitrary long, we say that $\conf{c}$ has infinite
level.
\end{defi}

It can be proved that $\espace$ contains infinite decreasing chains for
$\preceq$ and, as proved in \cite{villejournalcb}, countable SFTs may contain
such chains too. As a consequence, levels may not always be
finite nor well defined.
On the other hand, if any configuration $\conf{c}$ of a subshift is either at
level $n$ or below or there exists a configuration $\conf{c'}\preceq\conf{c}$
at level~$n$ then level $n+1$ is well defined by applying Zorn's lemma.
For example, level~$0$ is well defined since for every configuration $\conf{c}$
the subshift $\engendre{c}$ contains a minimal configuration.
For level~$1$, we have to prove that any decreasing chain in the set of
configurations of the subshift that are not at level~$0$ has a lower bound: This
set is not compact and the argument used for proving the existence of level~$0$
does not apply.
Nevertheless, in some cases, we can prove that there exists no infinite
decreasing chain for $\preceq$ down to level~$0$: This implies that any
decreasing chain in the set of configurations of the subshift that are not at
level~$0$ has a lower bound and then Zorn's lemma ensures us this set has a
minimal element; these minimal elements are precisely the configurations at
level~$1$.
A corollary of Proposition~\ref{prop:schmidtcaracsft} will be useful
for the study of these infinite decreasing chains:

\begin{prop}
\label{prop:decchainnonsft}
Let $\conf{x}_i$ be an infinite strictly decreasing chain (for the pre-order
$\preceq$) of configurations, then $\sshift{X}=\left\{\conf{x}\in \espace |
\forall i, \conf{x}\preceq \conf{x}_i \right\}=\bigcap_{i}\engendre{x_i}$ is a
subshift and not a subshift of finite type.
\end{prop}

\subsection{Cantor-Bendixson rank}

Given a topological space $\setofconf{X}$, its \emph{topological derivative}
is the set of configurations of $\setofconf{X}$ that are not isolated in
$\setofconf{X}$;  we denote it by $\setofconf{X'}$.
If $\sshift{X}$ is a subshift then $\sshift{X'}$ is also a
subshift since it is closed and shift-invariant.
However, if $\sshift{X}$ is SFT then its derivative is not necessarily SFT and
can actually be much more complex~\cite{villecbcount,villejournalcb}.
Note that a configuration $\conf{c}$ is isolated in $\sshift{X}$ if and only
if it is an equicontinuity point for the $\Z^d$ action by shifts if and only if
$\conf{c}$ is the only configuration of $\sshift{X}$ containing a given pattern
$\motif{P}$. In this latter case, we say that $\motif{P}$ \emph{isolates}
$\conf{c}$.

\begin{defi}
The topological derivation is defined for every ordinal $\lambda$ by induction:
\begin{itemize}
\item $\setofconf{X}^{(0)}=\setofconf{X}$
\item $\setofconf{X}^{(\alpha+1)}=(\setofconf{X}^{(\alpha)})'$
\item $\setofconf{X}^{(\lambda)}=\bigcap\limits_{\alpha<\lambda}\setofconf{X}^{(\alpha)}$
\end{itemize}
\end{defi}

\begin{figure}[htb]
\begin{center}
\opt{dessinetikz}{\input dessins/rangs}\opt{pdftikz}{\opt{bw}{\includegraphics{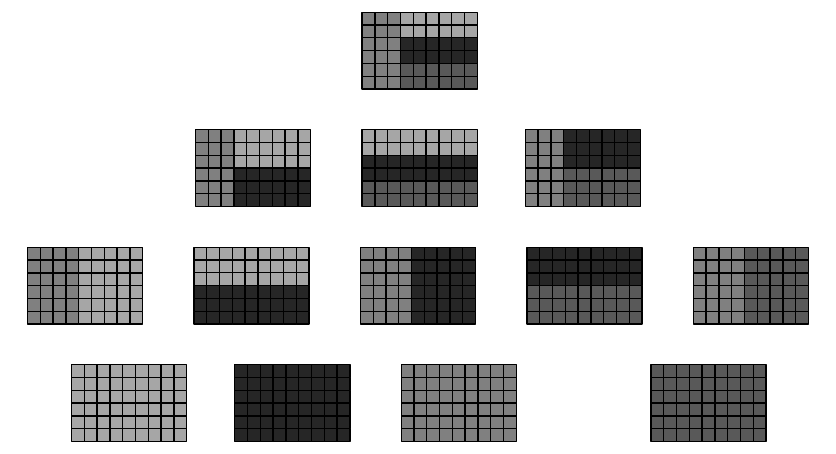}}\opt{color}{\includegraphics{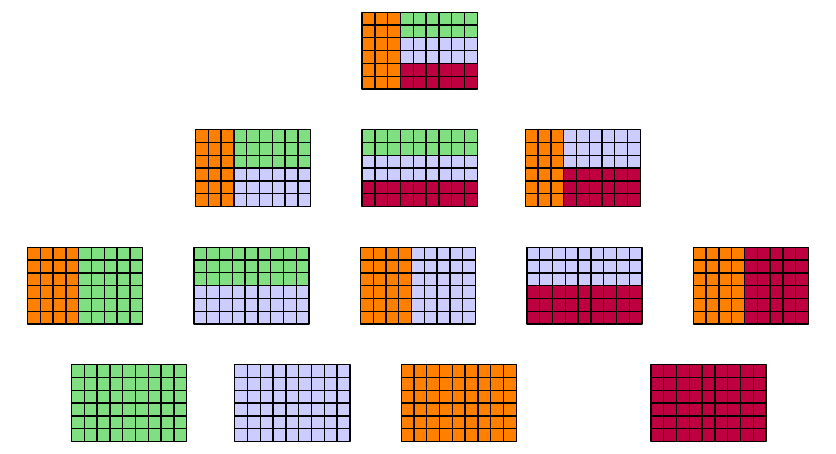}}}
\end{center}
\caption{Example of a sequence of topological derivations for an SFT.
The highest row represents isolated configurations, the next one those that are
isolated after removing the highest row, etc.
}
\label{fig:cbrangs}
\end{figure}

In the same way as for $\sshift{X'}$, if $\sshift{X}$ is a subshift then for any
ordinal $\lambda$, $\sshift{X}^{(\lambda)}$ is also a subshift since subshifts
are closed under intersection.
A sequence of topological derivations is represented on
Figure~\ref{fig:cbrangs}.

A configuration $\conf{c}$ of $\setofconf{X}$ has \emph{Cantor-Bendixson rank
$\beta$ in $\setofconf{X}$} if $\beta$ is the smallest ordinal such that
$\conf{c}\not\in \setofconf{X}^{(\beta)}$.
If there does not exist such a $\beta$, then we say that $\conf{c}$ has no 
Cantor-Bendixson rank.
We denote this rank $\rk{X}(\conf{c})$.

For any subshift $\sshift{X}$, there exists a \emph{countable} ordinal $\lambda$
such that $\sshift{X}^{(\lambda)}=\sshift{X}^{(\lambda+1)}$. Indeed, each step
of the topological derivation forbids patterns to appear and there is only a
countable number of patterns.
We call the smallest such ordinal the \emph{Cantor-Bendixson rank} of
$\setofconf{X}$~\cite[Paragraphe 24]{kura} and denote it by $\rk{}(\sshift{X})$.
As a consequence, \emph{if the Cantor-Bendixson rank is well defined for every
configuration of $\sshift{X}$ then $\sshift{X}$ is countable} as a countable
union of countable sets.
When this is not the case, $\sshift{X}^{(\rk{}(\sshift{X}))}$ is a non-empty
perfect subset of $\sshift{X}$ and, as such, has the cardinality of continuum. This
reasoning lead us to the following results:

\begin{prop}[\cite{stacsstructuralaspect}]
\label{prop:rankesssidnb}
The following are equivalent for a subshift $\sshift{X}$:
\begin{enumerate}
\item $\sshift{X}$ is countable;
\label{caracdnb:dnb}
\item There exists an ordinal $\lambda$ such that
$\sshift{X}^{(\lambda)}=\emptyset$;
\label{caracdnb:derivvide}
\item All the elements of $\sshift{X}$ have a Cantor-Bendixson rank.
\label{caracdnb:cbrank}
\end{enumerate}

As a consequence, a subshift is either finite, infinite countable or has the
cardinality of continuum.
\end{prop}

The purely topological operation of derivation of subshifts has
remarkably some intimate links with the combinatorial pre-order
$\preceq$:

\begin{prop}[\cite{stacsstructuralaspect}]
\label{prop:pluspetitalorsrkplusgrand}
Let $\conf{x}$ and $\conf{y}$ be two configurations being ranked by the
Cantor-Bendixson rank in $\sshift{X}$ such that $\conf{x}\prec \conf{y}$.
Then $\rk{X}(\conf{x})>\rk{X}(\conf{y})$.
\end{prop}

\begin{lemma}
\label{lemma:dersftdifffin}
If $\sshift{X}$ is a subshift and its topological derivative $\sshift{X}'$ is an
SFT then $\sshift{X}\setminus\sshift{X}'$ is finite up to translations.
\end{lemma}

\begin{proof}
Let $\left\{\motif{P}_n, n\in\N\right\}$ be the set of isolating patterns for
$\sshift{X}$. Let $\sshift{X}_n$ be the subshift where we forbid the patterns
forbidden by $\sshift{X}$ and, additionally, $\motif{P}_1,\ldots,\motif{P}_n$.
Clearly, $\cap_{n\in\N}\sshift{X}_n=\sshift{X}'$,
$\sshift{X}_{n+1}\subseteq\sshift{X}_n$ and since $\sshift{X}'$ is an SFT, by
Proposition~\ref{prop:schmidtcaracsft}, there exists $N\in\N$ such that
$\sshift{X}_N=\sshift{X}'$.

By definition of an isolating pattern, there is only one configuration of
$\sshift{X}$, up to translations, that contains a given isolating pattern.
Since any configuration of $\sshift{X}\setminus\sshift{X}'$ contains one of the
$\motif{P}_1,\ldots,\motif{P}_N$, $\sshift{X}\setminus\sshift{X}'$ is finite up
to translations.
\end{proof}

\section{Finite and uncountable subshifts}
\label{sec:finunc}
In this section we characterize finite and uncountable subshifts by means of the
configurations they contain.
The finite case was already dealt with in \cite{stacsstructuralaspect}:
\begin{theorem}[\cite{stacsstructuralaspect}]
\label{thm:finissitsper}
A subshift $\sshift{X}$ is finite if and only if it contains only periodic
configurations. It is actually an SFT.
\end{theorem}

Uncountable subshifts are very rich but we have seen that uncountability of a
subshift can only come from the existence of a non-empty perfect subshift.
We now show that non-empty perfect subshifts can only exist because of two
reasons: the existence of an infinite minimal subshift or of a \emph{complete
tree of patterns}.

\begin{defi}
For a subshift $\sshift{X}$, an infinite binary tree $\arbre{T}$
where each node is labelled by a pattern, is a \emph{complete tree of patterns
for $\sshift{X}$} if for every infinite path $(u_i)_{i\in\N}$ in $\arbre{T}$
there exists at least one configuration $\conf{c}$ in $\sshift{X}$ verifying the
following for every $u_i$ in the path:
\begin{itemize}
\item If the path goes through the left son of $u_i$ then the pattern labelled
by $u_i$ appears in $\conf{c}$
\item If the path goes through the right son of $u_i$ then the pattern labelled
by $u_i$ does not appear in $\conf{c}$.
\end{itemize}
\end{defi}

Having a complete tree of patterns for a subshift is a sufficient condition for
it to be uncountable: It gives us an explicit injection from $\cardcont$ into
$\sshift{X}$.
However, there exists uncountable SFTs admitting no complete tree of
patterns such as the $104$ tiles aperiodic one~\cite{Ollinger08}: It has
$\cardcont$ different configurations but all of them
contain exactly the same patterns which makes it impossible to have a complete
tree of patterns.

\begin{lemma}
\label{lemme:approxnonshiftealorsparfait}
If a subshift $\sshift{X}$ contains two configurations $\conf{x}$ and
$\conf{y}$ such that $\conf{x}\approx\conf{y}$ and $\conf{x}$ is not equal to
$\conf{y}$ modulo a shift then $\sshift{X}$ has the cardinality of continuum.
\end{lemma}

\begin{proof}
We prove that $\engendre{x}$ is a perfect set:
Suppose that $\engendre{x}$ contains an isolated point: $\conf{z}$.
There exists $v\in\plan$ such that $\conf{z}=\shift{v}{x}$ otherwise $\conf{z}$
would be limit of a (non-stationary) sequence of shifted of $\conf{x}$ and would
not be isolated.
Let $\motif{P}$ be a pattern isolating $\conf{z}$ in $\sshift{X}$, \ie such
that $\engendre{x}\cap\open{P}=\left\{\conf{z}\right\}$. $\motif{P}$ also
appears in $\conf{y}$, thus there exists $w\in\plan$ such that $\shift{w}{y}\in\open{P}$,
therefore $\shift{w}{y}=\conf{z}=\shift{v}{x}$ and $\conf{x}$ is a shifted of $\conf{y}$.
\end{proof}

The following corollary of Lemma~\ref{lemme:approxnonshiftealorsparfait} will be
useful later:
\begin{cor}
\label{cor:dnbapproxalorsshifte}
In a countable subshift, $\conf{x}$ and $\conf{y}$ contain the same patterns
($\conf{x}\approx\conf{y}$) if and only if $\conf{x}$ is a shifted of $\conf{y}$.
\end{cor}

\begin{theorem}
\label{thm:caracsshiftnondnb}
A subshift $\sshift{X}$ is uncountable if and only if one of the two following
conditions hold:

\begin{enumerate}
        \item There exists $\conf{x}$ and $\conf{y}$ in $\sshift{X}$ such that 
                $\conf{x}\approx\conf{y}$ and $\conf{x}$ is not a shifted of
                $\conf{y}$.
\label{caracsshiftnondnb:approxnonshifte}
\item There exists a complete tree of patterns for $\sshift{X}$.
\label{caracsshiftnondnb:arbrebinaire}
\end{enumerate}
\end{theorem}

\begin{proof}
        $\Leftarrow$: Point~\ref{caracsshiftnondnb:approxnonshifte} follows from 
        Lemma~\ref{lemme:approxnonshiftealorsparfait} and
        the definition of a complete tree of patterns gives an explicit
        injection from $\cardcont$ into $\sshift{X}$ for
        Point~\ref{caracsshiftnondnb:arbrebinaire}.

        $\Rightarrow$: Let $\sshift{X}$ be an uncountable subshift; it has a
        non-empty perfect subset; w.l.o.g, we can assume $\sshift{X}$ to be
        perfect.
Suppose point~\ref{caracsshiftnondnb:approxnonshifte} does not hold for
$\sshift{X}$; this means that the relation $\approx$ corresponds to the relation
``being equal up to a shift'' in $\sshift{X}$ (and thus in its subsets too).
For a given configuration, the set of its shifted is at most countably infinite
thus the equivalence classes for $\approx$ in $\sshift{X}$ and its subsets are
countable.

We will now do an induction on $\sshift{X}$; our invariant is that $\sshift{X}$
is a perfect set for which equivalence classes for $\approx$ are countable.

Consider the set $\setofconf{U}$ (possibly empty) of configurations of
$\sshift{X}$ that contain all the patterns of $\Lang{\sshift{X}}$. For a pattern
$\motif{P}$, let $\sshift{X}\setminus\{\motif{P}\}$ be the set of configurations
of $\sshift{X}$ not containing $\motif{P}$. By definition, we have:
\[
\sshift{X}=\setofconf{U}\cup\bigcup_{\motif{P}\in\Lang{\sshift{X}}}\sshift{X}\setminus\{\motif{P}\}
\]

$\setofconf{U}$ is countable since it corresponds to an equivalence class for
$\approx$.
Since $\sshift{X}$ is supposed to be uncountable (because it is perfect), there
exists $\motif{P}$ such that $\sshift{X}\setminus\{\motif{P}\}$ is uncountable
too. Therefore, $\sshift{X}\setminus\{\motif{P}\}$ has a non-empty perfect
subset. We continue the induction with that perfect
subset for the left son of the node we are constructing.

The set of configurations of $\sshift{X}$ that contain $\motif{P}$ contains the
open set $\open{P}$ of $\sshift{X}$ which is non-empty since we considered only
patterns appearing in at least one configuration of $\sshift{X}$.
Since $\open{P}$ is also closed, $\sshift{X}\cap\open{P}$ is a perfect set that
contains only configuration containing $\motif{P}$.
The induction continues with that perfect set for the right son of the node.

Iterating the induction defines a tree; each infinite path defines at least one
configuration, otherwise we would obtain an empty intersection of closed sets
which, by compactness, would be finite, \ie one of the non-empty perfect sets
we just constructed would be empty.
\end{proof}

\begin{figure}[htb]
\begin{center}
\opt{dessinetikz}{\input dessins/parabole}\opt{pdftikz}{\opt{bw}{\includegraphics{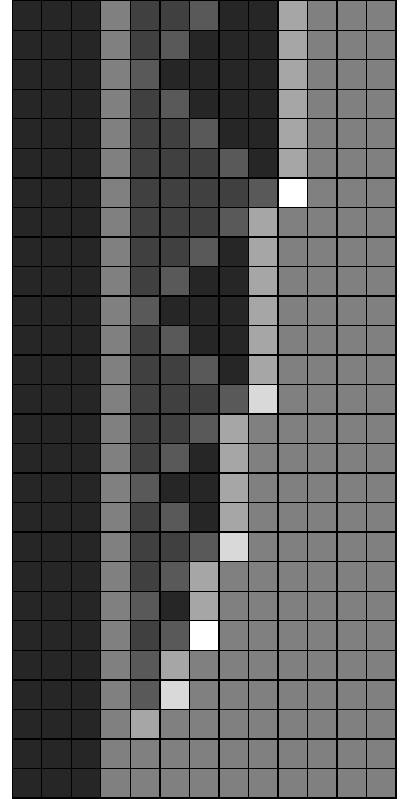}}\opt{color}{\includegraphics{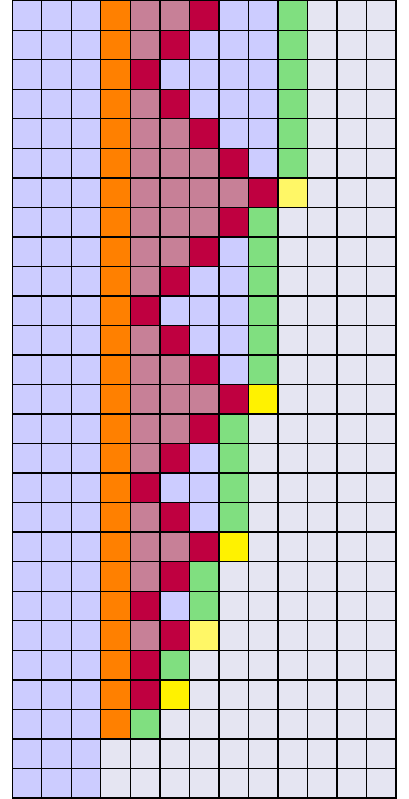}}}
\end{center}
\caption{Example of a typical configuration from a countable SFT verifying 
point~\ref{caracsshiftnondnb:arbrebinaire} but not
point~\ref{caracsshiftnondnb:approxnonshifte} of
Theorem~\ref{thm:caracsshiftnondnb}.}
\label{fig:parabole}
\end{figure}

We can construct uncountable SFTs verifying
point~\ref{caracsshiftnondnb:arbrebinaire} but not
point~\ref{caracsshiftnondnb:approxnonshifte}.
A vertical line ($\cellcouleur{\orange}$ on Figure~\ref{fig:parabole}) can appear
and forces a signal ($\cellcouleur{\purple}$ on Figure~\ref{fig:parabole}) that
will bump between the vertical line that serves as a sign post and another
signal that will draw a parabola ($\cellcouleur{\mgreen}$ on Figure~\ref{fig:parabole}).
When the signal hits the parabola, the latter shifts one cell to the right; the
signal has maximal speed which ensures us we are drawing a parabola.
Moreover, when the signal hits the parabola, the latter can use two different
colors when shifting to the right ($\cellcouleur{\yellow}$ and
$\cellcouleur{\myellow}$ cells  on Figure~\ref{fig:parabole}).
Since there are infinitely many cells where this may happen, we obtain an
infinite number of choices where to choose between two cells which ensures that
the obtained subshift has the cardinality of continuum.

Finally, if we know all the patterns of a configuration, since the number of
$\cellcouleur{\mgreen}$ cells above the $\cellcouleur{\yellow}$ and
$\cellcouleur{\myellow}$ ones allows us to determine uniquely the position of a
pattern, we can reconstruct (up to a shift) a configuration from its language.
Such a configuration cannot verify point~\ref{caracsshiftnondnb:approxnonshifte}
from Theorem~\ref{thm:caracsshiftnondnb} but the subshift is uncountable.

\section{The $\preceq$ pre-order for countable subshifts}
\label{sec:dnb}

Proposition~\ref{prop:rankesssidnb} states that a subshift is countable if and only
if all its elements have a Cantor-Bendixson rank.
This result, albeit very useful and used in this section, does not really help
us in understanding the structure of countable subshifts.
We study here in details this structure.

Since every subshift contains a minimal subshift, it also contains a minimal
element $\conf{x}$ for $\preceq$.
If $\conf{x}$ is not periodic, $\engendre{x}$ is infinite and thus perfect
(otherwise it would contain a non-empty subshift and thus not be minimal).
Thus, Lemma~\ref{lemme:approxnonshiftealorsparfait} ensures us that a minimal
subshift of a countable subshift contains only periodic configurations:

\begin{prop}[\cite{Dolbilin95,stacsstructuralaspect}]
\label{prop:minpourdnbper}
Every minimal element for $\preceq$ in a countable subshift is periodic.
Equivalently, every minimal subshift of a countable subshift is finite.
\end{prop}

Recall that an infinite subshift contains a configuration that is not periodic
by Theorem~\ref{thm:finissitsper}, thus an infinite countable subshift contains
a configuration that is not minimal (hence not at level~$0$) for $\preceq$ in
the subshift.
We now want to show that there is such a configuration at level~$1$, for which
we need to prove that the notion of level~$1$ is well defined. These results
were already proved in \cite{stacsstructuralaspect} but we include the proof
since it helps in understanding the approach used in this section.

\begin{lemma}[\cite{stacsstructuralaspect}]
\label{lemme:pasdechaineinfiniedecdanssfdnb}
In a countable subshift there exists no infinite decreasing chain for $\preceq$
to level~$0$.
\end{lemma}

\begin{proof}
Suppose that $(\conf{c_i})_{i\in\N}$ is such an infinite decreasing chain and
let $\sshift{X}=\bigcap_{i\in\N}\engendre{c_i}$.
$\sshift{X}$ is a subshift containing only periodic configurations, it is thus
a finite SFT by Theorem~\ref{thm:finissitsper}.
Therefore, by Proposition~\ref{prop:schmidtcaracsft}, there exists $N\in\N$ such
that $\sshift{X}=\engendre{c_N}$ and we supposed the sequence
$(\conf{c_i})_{i\in\N}$ to be strictly decreasing, a contradiction.
\end{proof}

\begin{lemma}[\cite{stacsstructuralaspect}]
\label{lemme:dnblvl1exists}
An infinite countable subshift $\sshift{X}$ contains a configuration at
level~$1$ for $\preceq$.
\end{lemma}

\begin{proof}
Let $\sshift{X}^{\clubsuit}$ be the set of non minimal configurations of
$\sshift{X}$ and let $(\conf{c_i})_{i\in I}$ be a decreasing chain in
$\sshift{X}^{\clubsuit}$. If it is infinite, by
Lemma~\ref{lemme:pasdechaineinfiniedecdanssfdnb}, the set of configurations
$\conf{c}\in\sshift{X}$ such that $\forall i\in I, \conf{c}\preceq\conf{c_i}$
cannot contain only minimal configurations, therefore, there exists
$\conf{c'}\in\sshift{X}^{\clubsuit}$ such that $\forall i\in I,
\conf{c'}\preceq\conf{c_i}$; if the chain is finite, it is straightforward that
such a $\conf{c'}$ exists.
As a consequence, every decreasing chain in $\sshift{X}^{\clubsuit}$ has a
lower bound.

By Zorn's lemma, $\sshift{X}^{\clubsuit}$ contains a minimal element for
$\preceq$. This minimal element is such that every configuration strictly
smaller for $\preceq$ is minimal in $\sshift{X}$, it is therefore a
configuration at level~$1$ for $\preceq$ in $\sshift{X}$.
\end{proof}

Since configurations at level~$1$ do exist in countable subshifts, we can
characterize them and obtain as a corollary the existence of these
configurations in a countable subshift:

\begin{theorem}
\label{thm:caraclvl1}
In a countable subshift $\sshift{X}$, a configuration $\conf{c}$ is at level~$1$
if and only if it is one of the following types:
\begin{enumerate}
\item all the patterns in $\conf{c}$ appear infinitely many times: in that case,
$\conf{c}$ admits one direction of periodicity (but not two) and is equal on
two half planes with one (one for each half plane) configuration
$\conf{x}\prec\conf{c}$
\label{caraclvl1:vp}
\item $\conf{c}$ contains a pattern that appears only once: in that case,
$\conf{c}$ is equal to a periodic configuration everywhere but on a finite
domain.
\label{caractlvl1:persffini}
\end{enumerate}
\end{theorem}

\begin{proof}
First, if $\conf{c}$ is a configuration matching either
condition~\ref{caraclvl1:vp} or condition~\ref{caractlvl1:persffini}, it is
straightforward that any converging sequence $(\shift{v_i}{c})_{i\in\N}$
converges either to a periodic configuration or to a shifted of $\conf{c}$,
hence $\conf{c}$ is at level~$1$.
Moreover, either every pattern in $\conf{c}$ appears infinitely many times or
$\conf{c}$ has a pattern that appears only once, hence we only have to prove
the conclusions of both cases.

By Proposition~\ref{prop:rankesssidnb}, there exists an ordinal $\lambda$ such that
$\conf{c}$ is isolated in $\derivee{X}{\lambda}$.
Let $\motif{P}$ be a pattern isolating $\conf{c}$ in $\derivee{X}{\lambda}$.
$\sshift{X}^{\clubsuit}=\derivee{X}{\lambda+1}\cap\engendre{c}$ is a subshift
containing exactly the configurations $\conf{x}\prec\conf{c}$, thus only periodic
configurations by Proposition~\ref{prop:minpourdnbper}.
$\sshift{X}^{\clubsuit}$ is a finite SFT by Theorem~\ref{thm:finissitsper}.
Let $\ForbPat$ be a finite set of pattern defining $\sshift{X}^{\clubsuit}$: $\Tilings =
\sshift{X}^{\clubsuit}$.

Case~\ref{caraclvl1:vp}:
All the patterns of $\conf{c}$ appear infinitely many times, therefore
$\motif{P}$ appears infinitely many times in $\conf{c}$,
it thus appears at the same position in $\shift{v}{c}$ for a $(0,0)\neq v\in\plan$.
Since $\motif{P}$ isolates $\conf{c}$ (and is also in $\derivee{X}{\lambda}$), 
$\shift{v}{c}=\conf{c}$. Hence, $\conf{c}$ admits one direction of
periodicity.
Now, if $\conf{c}$ admitted two directions of periodicity, $\conf{c}$ would be
periodic and thus at level~$0$.

Let $\conf{y}$ be a configuration of $\Tilings$; $\conf{y}$ is periodic and
$\conf{y}\prec\conf{c}$. Since we already proved that $\conf{c}$ admits a
direction of periodicity, we only have to prove that $\conf{c}$ is equal to
$\conf{y}$ on a quarter of plane.
For a given integer $n$, let $\motif{A_n}$ be the square pattern defined on
$[-n;n]^2$ in $\conf{y}$.
$\motif{A_n}$ appears somewhere in $\conf{c}$ since $\conf{y}\prec\conf{c}$.
Let $\motif{P_n}$ be a square pattern defined on $[-k_n;k_n]^2$ in $\conf{c}$
containing $\motif{A_n}$ and extended as much as possible
while still not containing a pattern of $\ForbPat$.
Remark that we cannot indefinitely extend it, and therefore $k_n$ is well
defined, otherwise we would have that
$\conf{c}\preceq\conf{y}$.

There exists a point at the border of $\motif{P_n}$ where a pattern of
$\ForbPat$ appears; let $\conf{c_n}$ be a shifted of $\conf{c}$ centered on this
pattern.
By compactness, w.l.o.g., we can assume that the sequence
$(\conf{c_n})_{n\in\N}$ tends to $\conf{z}$.
We can also assume that $\motif{P_n}$ always appears in the same quarter of
plane in $\conf{c_n}$.
Since the size of $\motif{P_n}$ tends to infinity when $n$ tends to infinity (by
definition $k_n \geq n$), $\conf{z}$ is equal to $\conf{y}$ on a quarter of
plane.
Moreover, each $\conf{c_n}$ contains at position $(0,0)$ a pattern of
$\ForbPat$; since the patterns of $\ForbPat$ are of bounded size, $\conf{z}$,
being a limit of $(\conf{c_n})_{n\in\N}$, also contains a pattern of $\ForbPat$
at position $(0,0)$.
We thus have $\conf{y}\prec\conf{z}\preceq\conf{c}$ and since we supposed
$\conf{c}$ to be at level~$1$, $\conf{z}\approx\conf{c}$.
By Corollary~\ref{cor:dnbapproxalorsshifte}, $\conf{z}$ is a shifted of 
$\conf{c}$ and thus $\conf{c}$ contains a quarter of plane that is equal to 
$\conf{y}$.
The direction of periodicity of $\conf{c}$ allows us to extend this quarter of
plane to a half plane.

On the other half plane, we can extract a sequence of shifted of
$\conf{c}$ that tends to a periodic configuration (since the limit does not
contain $\motif{P}$).
By the same reasoning, we conclude that $\conf{c}$ is equal to a periodic
configuration on half a plane disjoint to the first one.

Case~\ref{caractlvl1:persffini}:
There exists a pattern that appears only once in $\conf{c}$.
W.l.o.g., up to extending it, we can assume that this pattern is isolating
$\conf{c}$ and we denote it by $\motif{P}$.
Let $\ensdef{K}\subseteq\plan$ be the set of points where a pattern of
$\ForbPat$ appears in $\conf{c}$.
If $\ensdef{K}$ is infinite, there exists a sequence of shifted of
$\conf{c}$ that all contain a pattern of $\ForbPat$ at their center; by
compactness, w.l.o.g., this sequence converges towards a limit denoted by
$\conf{x}$. By construction, $\conf{x}$ is not in $\sshift{X}^{\clubsuit}$ since
it contains a pattern of $\ForbPat$ at its center; moreover, $\conf{x}$ does not
contain $\motif{P}$ and therefore $\conf{x}$ is a configuration that is not
minimal but such that $\conf{x}\prec\conf{c}$ which is impossible since we
supposed $\conf{c}$ at level~$1$.

Since $\sshift{X}^{\clubsuit}$ is finite, there exists a finite domain
$\ensdef{D}$ of $\plan$ such that if two configurations of
$\sshift{X}^{\clubsuit}$ match on $\ensdef{D}$ then they are equal.
Since $\ensdef{D}$ and $\ensdef{K}$ are finite, the set $\ensdef{E} =
\left\{a\in\plan, \exists b\in\ensdef{D}, a+b\in\ensdef{K}\right\}$ is also
finite; w.l.o.g., we can assume $\plan\setminus\ensdef{E}$ connected.
Since $\conf{c}$ contains no pattern of $\ForbPat$ on
$\plan\setminus\ensdef{E}$, $\conf{c}$ is therefore equal to a configuration of
$\sshift{X}^{\clubsuit}$ on $\plan\setminus\ensdef{E}$ which completes the
proof.
\end{proof}

\begin{cor}
\label{cor:dnbniveau1typeadoncauplusdeuxconfdessous}
In a countable subshift, configurations at level~$1$ have, up to shifts, at most
two (minimal) configurations below them for $\preceq$.
\end{cor}

Configurations at level~$1$ in a countable subshift are now characterized, the
only two possible configurations are depicted on Figure~\ref{fig:typelvl1}.
The fact that these configurations are equal on half planes to a periodic
configuration yields that they are never much more complicated that the simple
drawings on Figure~\ref{fig:typelvl1}.
Moreover, some cases never happen when we are restricted to SFTs instead of
subshifts:

\begin{figure}[htb]
\begin{center}
\opt{dessinetikz}{\input dessins/typelvl1}\opt{pdftikz}{\opt{bw}{\includegraphics{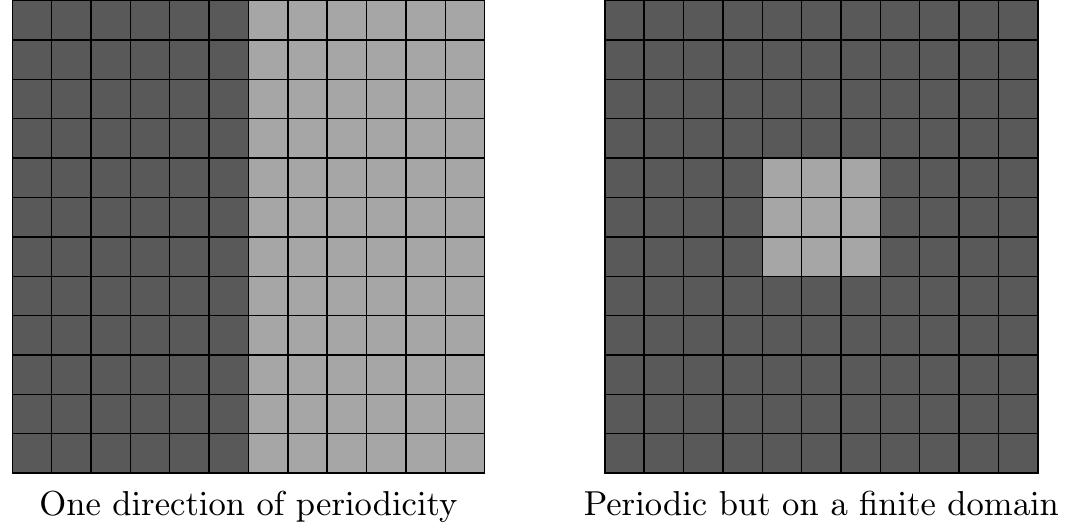}}\opt{color}{\includegraphics{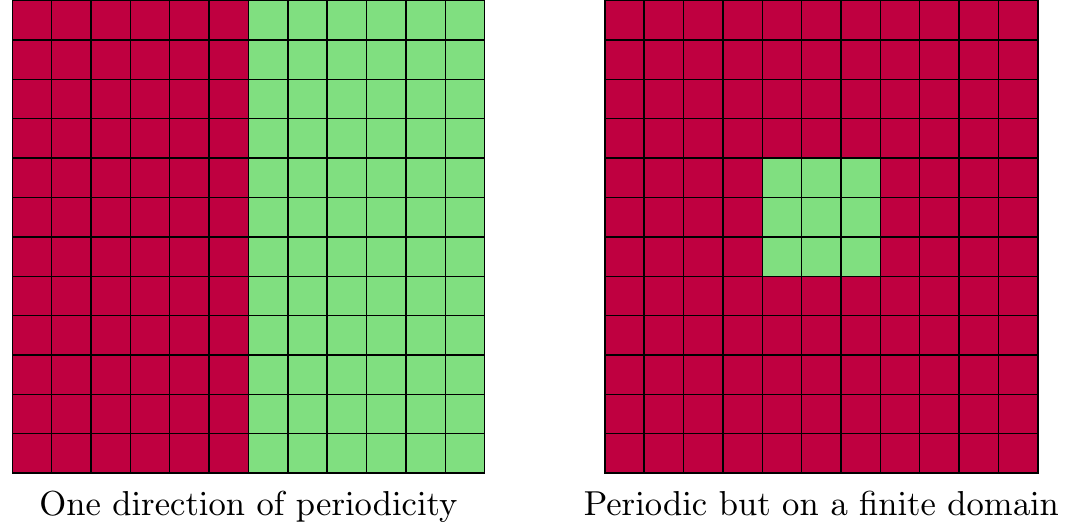}}}
\end{center}
\caption{The only two possible configurations at level~$1$ in a countable
subshift}
\label{fig:typelvl1}
\end{figure}

\begin{theorem}[\cite{stacsstructuralaspect}]
\label{thm:sftdnbalorsconf1per}
Every pattern in a configuration at level~$1$ in a countable SFT appears
infinitely many times.
Every infinite countable SFT contains a configuration that admits
exactly one direction of periodicity.
\end{theorem}

\begin{figure}[htb]
\begin{center}
\opt{dessinetikz}{\input dessins/hassenonsft}\opt{pdftikz}{\opt{bw}{\includegraphics{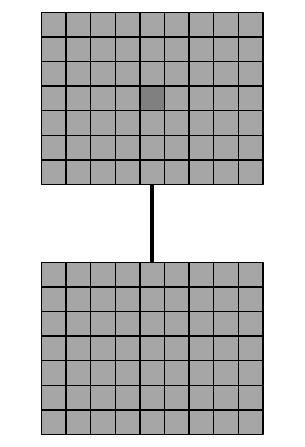}}\opt{color}{\includegraphics{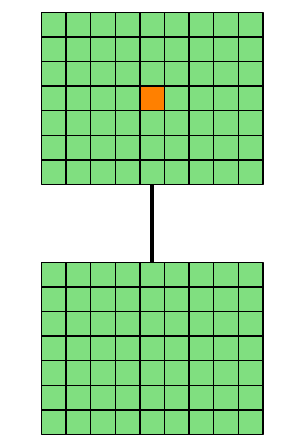}}}
\end{center}
\caption{Example of a Hasse diagram for $\preceq$ of a countable subshift not
        agreeing with Theorem~\ref{thm:sftdnbalorsconf1per}.
}
\label{fig:hassenonsft}
\end{figure}

Remark that Theorem~\ref{thm:sftdnbalorsconf1per} does not hold for general
subshifts:
Let $\alphab=\{
\tikz[scale=.2]{\filldraw[fill=\mgreen] (0,0) rectangle +(1,1);}
,
\tikz[scale=.2]{\filldraw[fill=\orange] (0,0) rectangle +(1,1);}
\}$. The set of forbidden patterns is the set of all patterns containing, at
least, two $\tikz[scale=.2]{\filldraw[fill=\orange] (0,0) rectangle +(1,1);}$
cells.
A Hasse diagram of the $\preceq$ pre-order, on which we can see that a
configuration at level~$1$ has a pattern appearing only once is depicted on
Figure~\ref{fig:hassenonsft}.

\begin{lemma}
\label{lemme:lvl1sftabc}
In a countable SFT, for a configuration $\conf{x}$ at level~$1$ with a
direction of periodicity $v$, any bi-infinite word $(\conf{x}(a+iw))_{i\in\Z}$
is of the form $^\omega\motif{A}\motif{B}\motif{C}^\omega$ ($a,w\in\Z^2$) as soon as
$w$ is not parallel to $v$.

Moreover, we can assume $\motif{A}, \motif{B}$ and $\motif{C}$ to be defined on
the same domain and $\left\{0,v\right\}\subseteq
\domain{A}=\domain{B}=\domain{C}$.
\end{lemma}

\begin{proof}
Let $\conf{c}\in\alphab^{\Z}$ be such a bi-infinite word.
By Theorem~\ref{thm:caraclvl1}, $\conf{x}$ contains two periodic half planes,
therefore $\conf{c}$ is ultimately periodic in both positive and negative
directions (since $v$ is not parallel to $w$) and is thus of the form
$^\omega\motif{A}\motif{B}\motif{C}^\omega$.

For having $\motif{A}, \motif{B}$ and $\motif{C}$ containing both $0$ and $v$
it suffices to consider the bi-infinite word in the first part of the proof to
have its letters in $\alphab^{D}$ where $D$ is a square domain containing
$0$ and $v$.
For having them defined on the same domain, it suffices to take a multiple of
the sizes of those three words.
\end{proof}

Remark that with Lemma~\ref{lemme:lvl1sftabc}, knowing $\motif{A}$, $\motif{B}$,
$\motif{C}$ and $v$ allows to uniquely reconstruct $\conf{x}$ and hence
configurations at level~$1$ in a countable SFT are computable.

\begin{cor}
\label{cor:sftdnbnv12dessous}
A configuration $\conf{c}$ at level~$1$ in a countable SFT has exactly two
configurations (modulo shift) below it.
\end{cor}

\begin{proof}
By Lemma~\ref{lemme:lvl1sftabc} we can write $\conf{c}$ as
$^\omega\motif{A}\motif{B}\motif{C}^\omega$. If $\motif{A}=\motif{C}$, then
$\motif{A}\neq\motif{B}$ otherwise $\conf{c}$ would be periodic and thus at
level~$0$. In this case, $(\motif{A}\motif{D}_i)_{i\in\Z}$ would also be a valid
configuration of our SFT with $\motif{D}_i\in\left\{\motif{A},
\motif{B}\right\}$, and our SFT would be uncountable.

Since $\motif{A}\neq\motif{C}$, these two patterns give us two minimal
configurations below $\conf{c}$ that are different modulo shift. By
Corollary~\ref{cor:dnbniveau1typeadoncauplusdeuxconfdessous}, they are the only
ones.
\end{proof}

\begin{figure}[htb]
\begin{center}
\opt{dessinetikz}{\input dessins/hassenonsft2}\opt{pdftikz}{\opt{bw}{\includegraphics{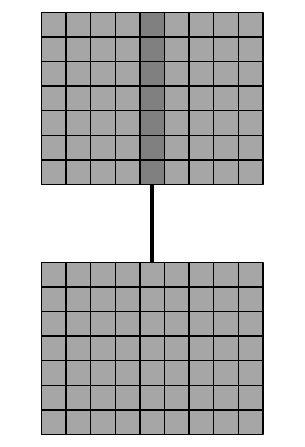}}\opt{color}{\includegraphics{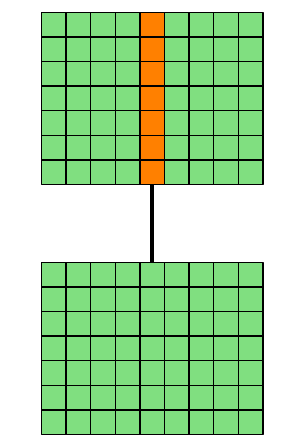}}}
\end{center}
\caption{Example of a Hasse diagram of a countable subshift.}
\label{fig:hassenonsft2}
\end{figure}

The subshift depicted on Figure~\ref{fig:hassenonsft2}
is an example of a countable (sofic) subshift with a configuration at level~$1$ in which every
pattern appears infinitely many times and with only one configuration below it
for $\preceq$, showing that Corollary~\ref{cor:sftdnbnv12dessous} cannot be
extended to arbitrary subshifts.

We proved that there exists no infinite decreasing chain down to level~$0$ in a
countable subshift. \cite[Theorem~$5.2$]{villejournalcb} gives a construction of
a countable SFT with an infinite decreasing chain down to level~$2$ and the
authors showed us (personal communication) how to obtain a countable sofic shift with an infinite
decreasing chain down to level~$1$. Therefore, the only remaining case is the
possibility to have an infinite decreasing chain down to level~$1$ in a
countable SFT. The following theorem proves that it is not possible:

\begin{theorem}
\label{thm:lvl2exists}
In a countable  SFT, there exists no infinite decreasing chain down to level~$1$
for $\preceq$.
\end{theorem}

In the proof of Theorem~\ref{thm:lvl2exists} we need a couple of lemmas; their
proof is provided after the theorem's proof for not burdening it more than it
needs to be.

\begin{proof}
Suppose $(\conf{c_n})_{n\in\N}$ to be such a chain.
Let $\sshift{S}=\bigcap_{n\in\N}\engendre{c_n}$ be the subshift that contains
(by definition) all the configurations that are below every $\conf{c_n}$ for
$\preceq$.
By hypothesis, $\sshift{S}$ contains only configurations at level~$0$ or $1$ and
as such only configurations with at least one direction of periodicity because
we are in an SFT.

Let $\conf{a}$ be a configuration at level~$1$ in $\sshift{S}$ and $v\neq (0,0)$
be its direction of periodicity.
Let $\motif{A}$, $\motif{B}$ and $\motif{C}$ be such as in
Lemma~\ref{lemme:lvl1sftabc} for $\conf{a}$ with, either an horizontal or
vertical vector since (at least) one of these two vectors is not parallel to
$v$.
Suppose $\conf{a}$ to be centered on $\motif{B}$ and let $\motif{M_k}$ be the
pattern defined on $[-k;k]^2$ at the center of $\conf{a}$.
By hypothesis, one can find $\motif{M_k}$ in every $\conf{c_n}$ since
$\conf{a}\preceq\conf{c_n}$.
We can obtain a more precise result on how to find these $\motif{M_k}$ in
$\conf{c_n}$:

\begin{lemma}
\label{lemme:bornesurlesmk}
For every $n$, there exists $N_n$ such that for any $\motif{M_{i}}$, $i\geq
N_n$, found in $\conf{c_n}$, the pattern $\motif{M_i}$ cannot be repeated
infinitely along $v$ and $-v$ to obtain an bi-infinite stripe.
\end{lemma}

In the following, we denote a bi-infinite periodic stripe of $\motif{M}$ along
$v$ and $-v$ by $\bandeinf{M}$; we say that a configuration $\conf{c}$ contains
$\bandeinf{M}$ if such an infinite stripe appears in $\conf{c}$.
In Lemma~\ref{lemme:bornesurlesmk} we can assume $N_n$ to be minimal, that is,
such that $\bandeinf{(\motif{M_{N_n-1}})}$ appears in $\conf{c_n}$.

\begin{lemma}
\label{lemme:kbornemiavecmkaucentrenonext}
There exists $k$ such that, in every $\conf{c_n}$, one can $\motif{M_i}$'s for
arbitrary large $i$'s such that the pattern $\motif{M_k}$ at the center of
$\motif{M_i}$ is not part of a $\bandeinf{M_k}$.
\end{lemma}

Let $k$ be such as in Lemma~\ref{lemme:kbornemiavecmkaucentrenonext}.
Let $\motif{M_i}$, $i\geq n$, be a pattern of $\conf{c_n}$ such that the pattern
$\motif{M_k}$ at its center cannot be extended to a $\bandeinf{M_k}$.
Let $\conf{x_n}$ be a shifted of $\conf{c_n}$ centered on the pattern
$\motif{M_k}$ that cannot be extended to a $\bandeinf{M_k}$.
By compactness one can obtain a ``limit'' of these $\conf{x_n}$: Since $k$ is
fixed and the $\motif{M_i}$'s are of strictly increasing sizes, we obtain a
limit $\conf{x}$ with a semi-infinite stripe consisting of a pattern
$\motif{A}^l\motif{B}\motif{C}^l$ repeated along, say, $v$, as depicted on
Figure~\ref{fig:lvl2existsbandesemiinf}.
Moreover, this configuration $\conf{x}$ admits a direction of periodicity since
it is a configuration at most at level~$1$ in a countable SFT.

\begin{figure}[htb]
\begin{center}
\opt{dessinetikz}{\input dessins/lvl2existsbandesemiinf}\opt{pdftikz}{\opt{bw}{\includegraphics{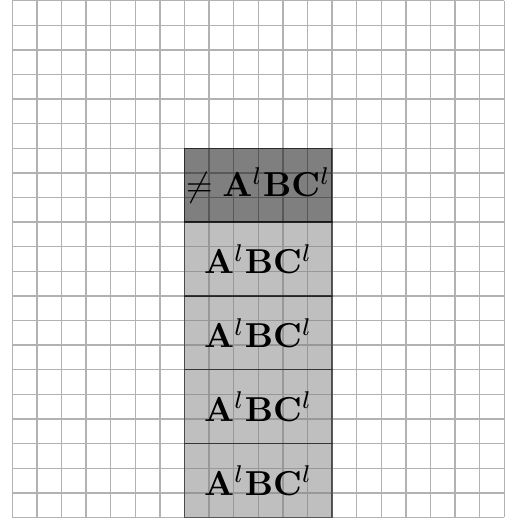}}\opt{color}{\includegraphics{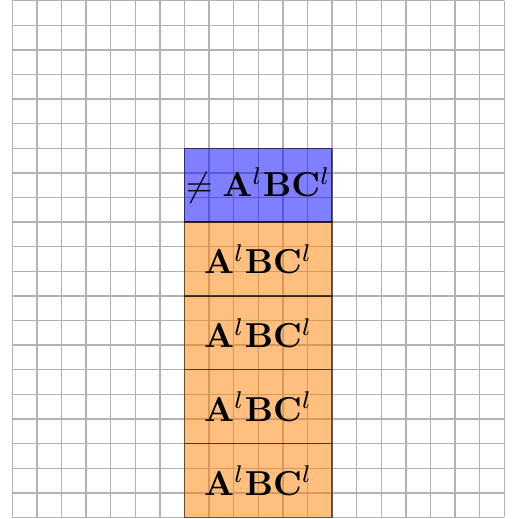}}}
\end{center}
\caption{Example of configuration extracted from $(\conf{x_{n}})_{n\in\N}$.}
\label{fig:lvl2existsbandesemiinf}
\end{figure}

The periodicity vector of $\conf{x}$ cannot be parallel to $v$ otherwise, by
periodicity, it would contain a $\bandeinf{M_k}$.
Therefore, up to multiplying the periodicity vector of $\conf{x}$ and taking a
``limit'' we can obtain an infinite stripe of type
$\bandeinf{(\motif{C}^l\motif{D}\motif{A}^l)}$.
Figure~\ref{fig:lvl2existsbandesemiinfvp} depicts $\conf{x}$ with a periodicity
vector $w$ and Figure~\ref{fig:lvl2existsbandesemiinfvpextr} the configuration
that can be extracted from $\conf{x}$.

\begin{figure}[htb]
\begin{center}
\opt{dessinetikz}{\input dessins/lvl2existsbandesemiinfvp}\opt{pdftikz}{\opt{bw}{\includegraphics{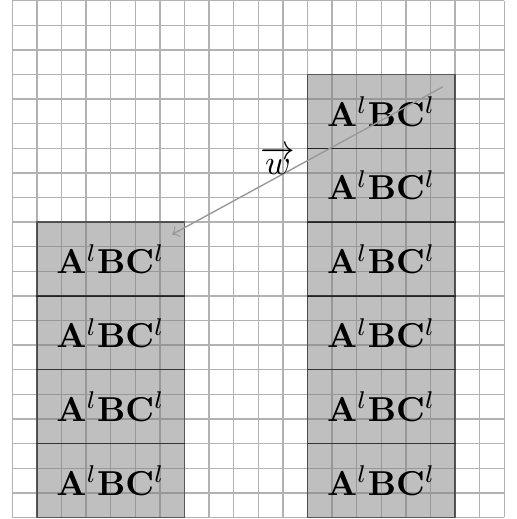}}\opt{color}{\includegraphics{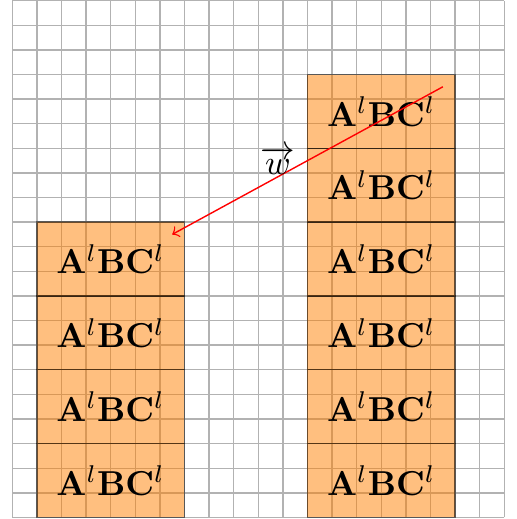}}}
\end{center}
\caption{Example of $\conf{x}$ and its periodicity vector $w$.}
\label{fig:lvl2existsbandesemiinfvp}
\end{figure}

\begin{figure}[htb]
\begin{center}
\opt{dessinetikz}{\input dessins/lvl2existsbandesemiinfvpextr}\opt{pdftikz}{\opt{bw}{\includegraphics{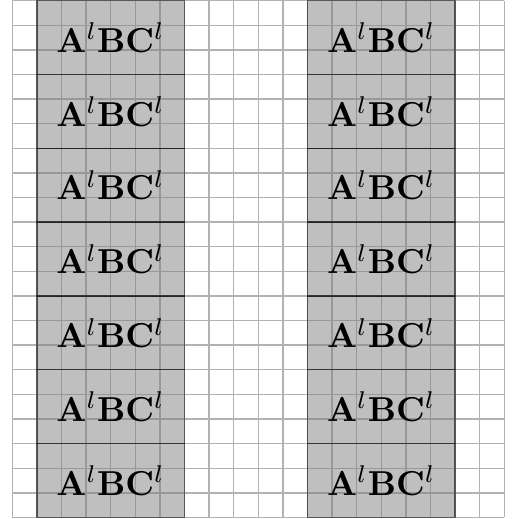}}\opt{color}{\includegraphics{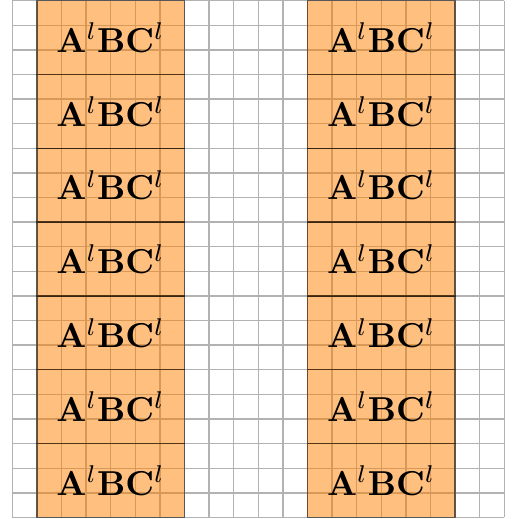}}}
\end{center}
\caption{Example of configuration that can be extracted from $\conf{x}$.}
\label{fig:lvl2existsbandesemiinfvpextr}
\end{figure}

If this were possible, we would be able to obtain $\cardcont$ different
configurations from $\conf{a}$.
This completes the proof: There exists no infinite decreasing chain for
$\preceq$ down to level~$1$ in a countable SFT.
\end{proof}

\begin{proof}[Proof of Lemma~\ref{lemme:bornesurlesmk}]
Suppose one can find $\bandeinf{M_i}$, for arbitrary larges $i$'s in
$\conf{c_n}$.
Fix such an $i$ with $\bandeinf{M_i}$ appearing in $\conf{c_n}$; since 
$\conf{c_n}$ is different from $\conf{a}$, there exits $j>i$ such that 
$\bandeinf{M_j}$ appears in $\conf{c_n}$ somewhere else from where 
$\bandeinf{M_i}$ appears.

We, therefore, obtain a pattern of type
$\motif{A}^l\motif{B}\motif{C}^l\motif{D}\motif{A}^k\motif{B}\motif{C}^k$ in
$\conf{c_n}$ (for example: $\motif{M_i}=\motif{A}^l\motif{B}\motif{C}^l$,
$\motif{M_j}=\motif{A}^k\motif{B}\motif{C}^k$ and $\motif{D}$ the separation
between them).
Thus, $\motif{C}^l\motif{D}\motif{A}^k\motif{B}\motif{C}^k$ is a valid pattern
for our SFT and $\conf{a}$ contains $\motif{C}^\omega$, so we have an infinite
number of positions where to replace a given $\motif{C}^m$ by
$\motif{C}^l\motif{D}\motif{A}^k\motif{B}\motif{C}^k$ to obtain $\cardcont$
different configurations in our SFT supposed to be countable.
\end{proof}

\begin{proof}[Proof of Lemma~\ref{lemme:kbornemiavecmkaucentrenonext}]
        First, if the sequence of minimal $(N_n)_{n\in\N}$ from
        Lemma~\ref{lemme:bornesurlesmk} is bounded, it suffices to take $k$
greater than all the $N_n$'s; we can therefore assume that $(N_n)_{n\in\N}$ can
be arbitrary large.

There exists a $k$, that does not depend
on $n$, such that all $\conf{c_n}$ contain at most one $\bandeinf{M_k}$: If two
existed, we would be able to find in $\conf{c_n}$ a pattern like 
$\motif{A}^i\motif{B}\motif{C}^i\motif{D}\motif{A}^j\motif{B}\motif{C}^j$ and
therefore obtain $\cardcont$ different configurations from $\conf{a}$ since we
are in an SFT.

Suppose we can find $\motif{M_i}$ in $\conf{c_n}$ for arbitrary large $i$'s along
$\bandeinf{M_k}$.
This stripe of $\motif{M_k}$'s is of bounded width by
Lemma~\ref{lemme:bornesurlesmk} and w.l.o.g. we can assume that the cell that
makes it impossible to extend indefinitely its width is always on the same side
of $\motif{A}$ in $\motif{M_{k}}$.
As such, if we have $i$ greater than the bound on the width of these stripes of 
$\motif{M_k}$ ($i\geq N_n$) we obtain a configuration where we have an important
number of $\motif{A}$ repeated along $v$ with a pattern different from
$\motif{A}$, denoted by $\motif{D}$, at the end. All of that is adjacent to a
$\bandeinf{M_{N_n-1}}$ as depicted on Figure~\ref{fig:lvl2existsmiailleurs}.

\begin{figure}[htb]
\begin{center}
\opt{dessinetikz}{\input dessins/lvl2existsmiailleurs}\opt{pdftikz}{\opt{bw}{\includegraphics{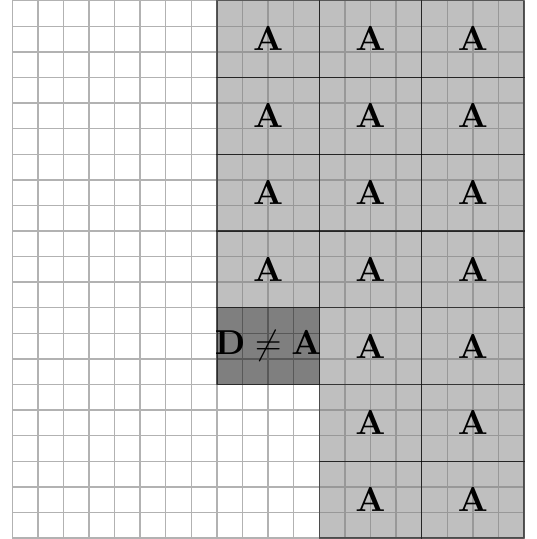}}\opt{color}{\includegraphics{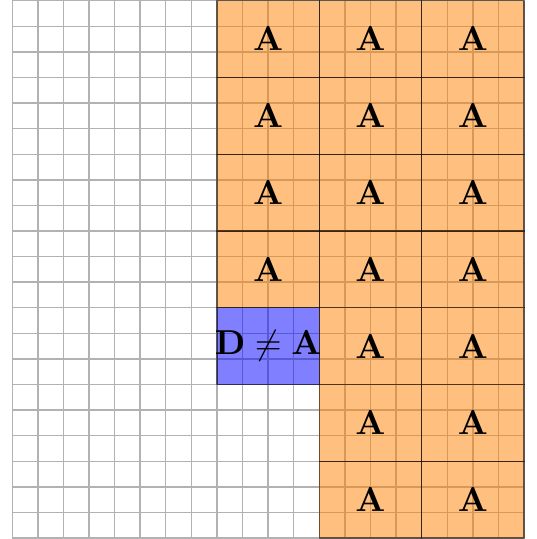}}}
\end{center}
\caption{Example of $\conf{c_{n}}$ when one can find arbitrary large
$\motif{M_i}$'s along a $\bandeinf{M_k}$.}
\label{fig:lvl2existsmiailleurs}
\end{figure}

Now, we can apply this reasoning to every $\conf{c_n}$ with patterns
$\motif{M_i}$ arbitrary large and since $N_n$ can be arbitrary large too we
obtain by taking a limit of $\motif{D}$ a configuration such as the one depicted
on Figure~\ref{fig:lvl2existsmiailleurs} but with patterns $\motif{A}$ repeated
infinitely.

This configuration admits a vector of periodicity since it is extracted from 
$\conf{c_n}$ but this periodicity would force to have 
$\motif{D}=\motif{A}$ which we supposed to be different.
\end{proof}

\begin{theorem}
\label{thm:lvl2ssi0vp}
A countable SFT contains a configuration without any direction of periodicity
if and only if it contains one at level~$2$ for $\preceq$.
\end{theorem}

\begin{proof}
$\Leftarrow:$
Theorem~\ref{thm:caraclvl1} ensures that a configuration at level~$2$ in a
countable subshift has no direction of periodicity.

$\Rightarrow:$ 
Let $\sshift{X}^{\clubsuit}$ be the set of configurations at level greater than
$1$ in the countable SFT.
It is non-empty since it contains all the configurations without any direction
of periodicity.
Theorem~\ref{thm:lvl2exists} states that there exists no infinitely decreasing
chain down to level~$1$ and Lemma~\ref{lemme:pasdechaineinfiniedecdanssfdnb} no
down to level~$0$. Therefore, any decreasing chain in $\sshift{X}^{\clubsuit}$
admits a lower bound.
As a consequence, by Zorn's lemma, there exists a minimal element in
$\sshift{X}^{\clubsuit}$: A configuration at level~$2$.
\end{proof}

\section{Cantor-Bendixson rank of countable subshifts}

When using the Cantor-Bendixson derivation for studying countable subshifts,
natural questions arise: What are the possible Cantor-Bendixson ranks of
countable SFTs of dimension~$2$? Can a countable SFT contain non computable
configurations? If yes, what is the minimal Cantor-Bendixson rank of such an
SFT?
At the time of writing \cite{stacsstructuralaspect} we did not even know if
there could exist countable SFTs with an infinite Cantor-Bendixson rank or with
non-computable configurations.
Shortly after, Ronnie Pavlov showed us an example of a countable sofic shift with
infinite rank. Later came the constructions from \cite{vanier2011} proving
that there exist two-dimensional SFTs with non-computable configurations and of
arbitrary large Cantor-Bendixson rank: 
The Cantor-Bendixson rank of an SFT is always a computable ordinal and a
consequence of \cite{vanier2011} is that for any computable ordinal there exists
an SFT of greater Cantor-Bendixson rank.
The Cantor-Bendixson rank of a countable effectively closed subset of
$\left\{0,1\right\}^{\N}$ can be any non-limit computable ordinal~\cite{pi01cbrank2,pi01cbrank}
and it has been recently proved that for any countable effectively
closed subset of $\left\{0,1\right\}^{\N}$ of rank $\lambda$, there exists a
countable SFT of rank $\lambda{}+4$~\cite[Theorem~$4.5$]{villejournalcb}
( a similar result but with a bigger increase in the Cantor-Bendixson rank can
also be found in \cite{vanier2011}).
Since it is easy to construct countable SFTs of rank $1,2,3,4$ and $5$ and a simple
compactness argument shows that the rank cannot be a limit ordinal, the only
ranks we do not know to exist are $\lambda+1$, $\lambda+2$, $\lambda+3$ and
$\lambda+4$ for $\lambda$ a limit ordinal.
In \cite{stacsstructuralaspect} we proved that the Cantor-Bendixson rank of a
countable subshift cannot be $\lambda+1$. We prove here that it cannot be $\lambda+2$
either, leaving open only the $\lambda+3$ and $\lambda+4$ cases.
We first use the Cantor-Bendixson derivatives to obtain the structure of the
configurations of interest here:

\begin{theorem}
\label{thm:rang2}
Let $\sshift{S}$ be a countable SFT and $\sshift{X}$ a subshift of $\sshift{S}$
of Cantor-Bendixson rank at most $3$, \ie $\sshift{X}^{(3)}=\emptyset$.
Then any configuration $\conf{c}$ of $\sshift{X}$ is of one of the following
types:
\begin{itemize}
\item $\conf{c}$ admits a direction of periodicity
\item $\conf{c}$ consists in a finite number of half-lines or full-lines and
none of them is parallel to another one. In this case, on every connected
component of $\plan$ that does not intersect with any of these lines,
$\conf{c}$ is equal to a bi-periodic configuration on this domain. 
Remark that since the number of lines is finite and there are no two parallel
lines, there is only a finite number of maximal (for the inclusion) such
connected components.
\end{itemize}
\end{theorem}

\begin{proof}
If $\sshift{X}^{(2)}=\emptyset$ then all the configurations of $\sshift{X}$ are
at most at level~$1$ by Proposition~\ref{prop:pluspetitalorsrkplusgrand} and
thus admit a direction of periodicity by Theorem~\ref{thm:sftdnbalorsconf1per}.
We can therefore assume that $\sshift{X}^{(2)}\neq\emptyset$, which means that
$\sshift{X}^{(2)}$ contains only periodic configurations and $\sshift{X}^{(1)}$
contains only configurations with one direction of periodicity.
$\sshift{X}^{(2)}$ is finite (because $\sshift{X}^{(3)}=\emptyset$), it is
therefore an SFT by Theorem~\ref{thm:finissitsper}.
$\sshift{X}^{(1)}$ is finite up to translations: $\sshift{X}^{(2)}$ is finite
and $\sshift{X}^{(1)}\setminus\sshift{X}^{(2)}$ is finite up to translations by
Lemma~\ref{lemma:dersftdifffin}.

Let $\conf{c}$ be a configuration in $\sshift{X}\setminus\sshift{X}^{(1)}$
without direction of periodicity (if there is no such configuration, our theorem
is proved).
$\conf{c}$ contains a pattern $\motif{M}$ that appears only once: $\conf{c}$ is
isolated in $\sshift{X}$ and were the pattern isolating it appearing twice,
$\conf{c}$ would have a direction of periodicity.

Since $\sshift{X}^{(1)}$ contains only configurations with one direction of
periodicity and is finite up to translations, there exists an integer $n$
greater than $3$ times the size of any vector of periodicity of any
configuration in $\sshift{X}^{(1)}$.
Let $E$ be the set of patterns defined on $[-n;n]^2$ appearing in a
configuration of $\sshift{X}^{(1)}$.

Let $N\subseteq\plan$ be the positions where a pattern that is not in $E$ appears in
$\conf{c}$. $N$ is finite: If it were infinite we could find a converging
sequence of shifted of $\conf{c}$ whose limit is different from $\conf{c}$
because it would not contain $\motif{M}$ but would contain a pattern that is not
in $E$; this limit would therefore not be isolated in $\sshift{X}$ but would not
belong to $\sshift{X}^{(1)}$ either.

Let $\motif{A}$ be a pattern defined on $[-n;n]^2$ appearing, outside of $N$, in
$\conf{c}$.
If $\motif{A}$ admits two different directions of periodicity then it appears in
a bi-periodic configuration of $\sshift{X}$ and we do not need to prove anything
for our theorem.
$\motif{A}$ admits a unique direction of periodicity, and a vector $v\neq(0,0)$ as a
period. In addition, suppose $\norme{}{v}$ to be minimal.
Let $\motif{B}$ be the pattern defined on $[-n;n]^2$ appearing in $\conf{c}$ at
the position shifted by $v$ from $\domain{A}$.
For now, assume that $x+v$ is not in $N$ either. $\motif{B}$ admits a vector $w$
as a period.

Let $x$ be a point in $\domain{A}$ such that $x+w$ and $x+v+w$ are also in
$\domain{A}$ (such a point exists because, by definition, $n\geq3\norme{}{w}$).
$x+v$ and $x+v+w$ are in $\domain{B}$, therefore
$\conf{c}(x+v+w)=\conf{c}(x+v)$. Since $\motif{A}$ admits $v$ as a period, we
obtain that $\conf{c}(x)=\conf{c}(x+w)$.
Since $x$ can be chosen arbitrary among the points in $[-n/2;n/2]^2$, and $n/2$
is greater than the size of any possible vector of periodicity, $\motif{A}$
admits also $w$ as a period.
Since $\motif{A}$ admits a unique direction of periodicity, there exists an
integer $y\geq1$ such that $w=yv$.

Let $x$ be in $\domain{B}$ such that $x+v$ is also in $\domain{B}$.
We get two cases:
\begin{itemize}
\item $x+v$ is in $\domain{A}$ (thus $x$ too): $\conf{c}(x)=\conf{c}(x+v)$ by
periodicity of $\motif{A}$.
\item $x+v$ is not in $\domain{A}$: $x+v-w$ and $x-w$ are in
$\domain{A}\cap\domain{B}$ since $w=yv$, then $\conf{c}(x) = \conf{c}(x-w) =
\conf{c}(x+v-w) = \conf{c}(x+v)$.
\end{itemize}

In any case, $\motif{B}$ also admits $v$ as a vector of periodicity. Therefore,
$\motif{A}$ and $\motif{B}$ are equal.
By an induction, this gives us a line if we never hit $N$ or an half-line if it
hits $N$ on one side.

Since $N$ is finite, there is a finite number of such half-lines that do not
contain patterns appearing in a bi-periodic configuration.
If a given (full-)line appears $3$ times in $\conf{c}$, on the same side of $N$,
then since $\sshift{S}$ is an SFT, we can find an infinite number of points
where to put or not this line, separated by the same bi-periodic patterns, and
thus $\sshift{S}$ would not be countable.

We can now conclude that the number of lines and half-lines is finite, the rest
being patterns appearing in a bi-periodic configuration.
Since we chose $n$ to be greater than any possible period, a straightforward
induction shows that $\conf{c}$ is equal to a bi-periodic configuration on any
connected component of $\plan$ that does not intersect any of these lines.
Since the number of lines is finite, the number of parallel lines is \latin{a
fortiori} finite too; we can therefore take the $n$ from this proof to be
sufficiently large so that the parallel lines are merged into one and get our
result.
\end{proof}

\begin{cor}
        \label{cor:compptrank4}
An SFT $\sshift{X}$ of Cantor-Bendixson rank at most $4$ contains only
computable configurations.
\end{cor}

The article \cite{pi01cbrank} proves that there exist countable effectively closed subset of
$\left\{0,1\right\}^{\N}$ with Cantor-Bendixson rank~$2$ containing
non-computable points; when combined with \cite[Theorem~$4.5$]{villejournalcb},
we obtain countable SFTs of Cantor-Bendixson rank~$6$ with non-computable
configurations hence after Corollary~\ref{cor:compptrank4} this leaves only
countable SFTs of rank~$5$ for which we do not know if they can contain non
computable configurations.

\begin{proof}
First, the isolated points of $\sshift{X}$ are computable: knowing the pattern
isolating them, it is possible to compute them since $\sshift{X}$ is compact
\cite{CenzerRec}.

Let $\sshift{X}'$ be the topological derivative of $\sshift{X}$. $\sshift{X}'$
matches the conditions of Theorem~\ref{thm:rang2}.
All the configurations described by the conclusions of Theorem~\ref{thm:rang2}
can be described by a finite amount of information and are therefore computable:
if the configuration does not admit any direction of periodicity it suffices to
know pattern defined over $N$, the position and the colors of the finite number
of half-lines and the colors of the bi-periodic parts separating them.
If the configuration admits exactly one direction of periodicity then, since it
is at level~$1$ by Theorem~\ref{thm:caraclvl1}, it can be described
by a finite amount of information by Lemma~\ref{lemme:lvl1sftabc}. If the
configuration is bi-periodic then knowing the pattern that is repeated
bi-periodically is sufficient in order to compute the whole configuration.
\end{proof}

\begin{figure}[htb]
\begin{center}
\opt{dessinetikz}{\input dessins/typic_rang_2_f}\opt{pdftikz}{\opt{bw}{\includegraphics{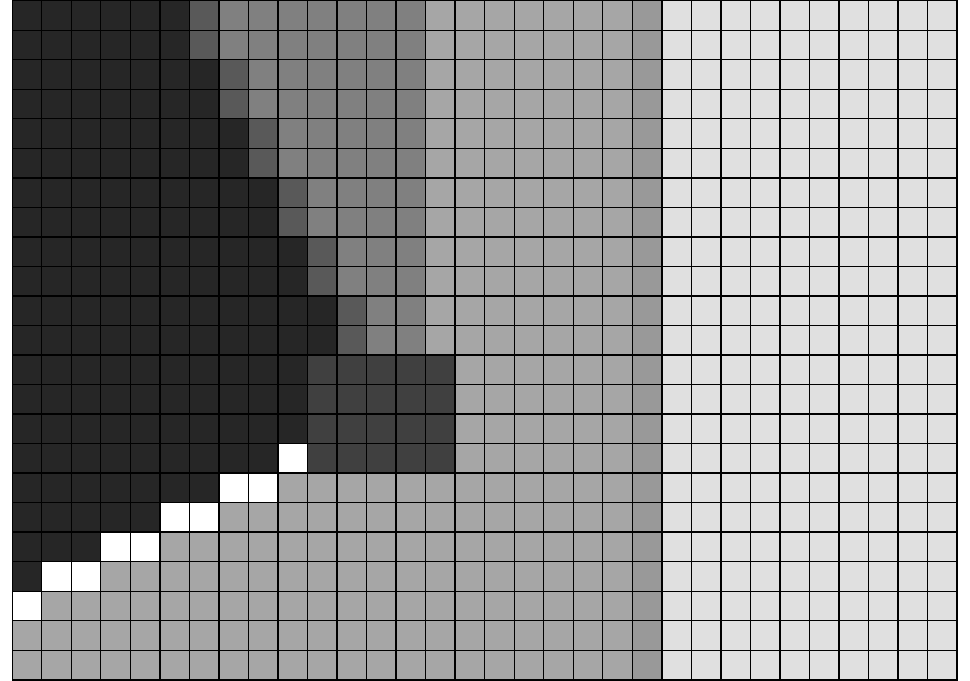}}\opt{color}{\includegraphics{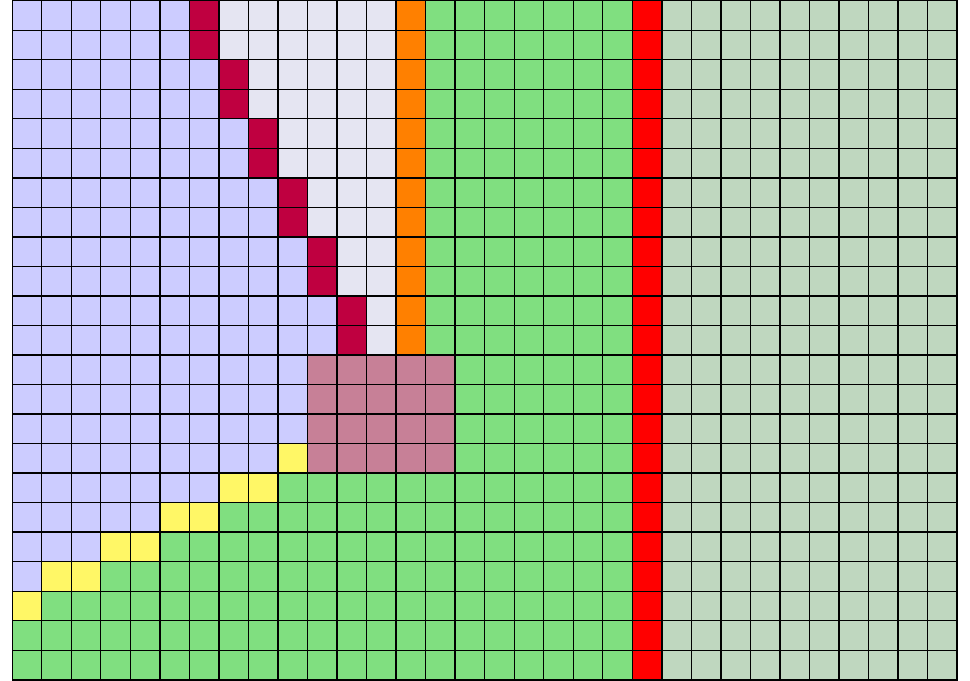}}}
\end{center}
\caption{A typical non-periodic configuration for Theorem~\ref{thm:rang2}.}
\label{fig:typic_rang_2}
\end{figure}

Theorem~\ref{thm:rang2} can be seen as the characterization of the
simplest, in the sense of the Cantor-Bendixson rank, aperiodic configurations in
a countable SFT; such a typical aperiodic configuration is depicted on
Figure~\ref{fig:typic_rang_2}.
As we have seen thorough this paper, the pre-order $\preceq$ shares intimate
links with the Cantor-Bendixson derivation and even stronger links when we are
dealing with countable SFTs. Nevertheless, even if we have been able to prove
that level~$2$ is always well defined in countable SFTs, we do not know how to
characterize such configurations at level~$2$. The natural conjecture would be
that those are precisely the non-periodic configurations matching the
conclusions of Theorem~\ref{thm:rang2}:

\begin{op}
\label{op:caraclvl2}
Characterize configurations at level~$2$ for $\preceq$ in countable SFTs.
Are they exactly the non-periodic configurations matching the conclusions of
Theorem~\ref{thm:rang2}?
\end{op}

We shall now prove that no countable SFT can be of Cantor-Bendixson rank
$\lambda+2$ for $\lambda$ a limit ordinal; for proving it, we first need a
lemma that might be of independent interest:

\begin{lemma}
    \label{lemma:perfinsft}
    Let $\mathcal{P}$ be a finite subset of $\plan$ 
    and $\sshift{X}$ be an SFT, then the set
    $\sshift{X}_{\mathcal{P}}\subseteq\sshift{X}$ of configurations of
    $\sshift{X}$ admitting a period in $\mathcal{P}$ is an SFT.
\end{lemma}

\begin{proof}
    It is clear that if $\mathcal{P}$ consists of only one vector then
    $\sshift{X}_{\mathcal{P}}$ is an SFT. We make an induction on $\mathcal{P}$
    by successive applications of the following claim:
\begin{claim}
    \label{claim:perunoncolpersft}
    Let $\mathcal{P}$ be a finite subset of $\plan$.
    Let $\sshift{X}$ be an SFT that contains only configurations admitting a
    period from $\mathcal{P}$. Let $q\in\plan$ be a vector that is not collinear
    to any vector from $\mathcal{P}$ and $\sshift{Y}_q$ be an SFT containing
    only configurations that are $q$-periodic.
    Then, $\sshift{X}\cup\sshift{Y}_q$ is an SFT.
\end{claim}

\begin{proof}[Proof of Claim~\ref{claim:perunoncolpersft}]
        Let $n\in\N$ be a fixed large integer; such that all the forbidden
        patterns defining $\sshift{X}$ or $\sshift{Y}_q$ are defined on
        $[-n/4;n/4]^2$ and all the vectors of 
        $\mathcal{P}\cup\left\{q\right\}$ are contained in $[-n/4;n/4]^2$.
        By a compactness argument, one can show that there exists $N$ such that
        any pattern defined on $[-N;N]^2$ not containing any
        $\sshift{X}-$forbidden pattern contains, at its center, a pattern
        defined on $[-n;n]^2$ appearing in a configuration of $\sshift{X}$: This
        pattern is $p-$periodic for some $p\in\mathcal{P}$.
        By the same argument, replacing $\sshift{X}$ with $\sshift{Y}_q$, we
        find a (maybe bigger) $N$ such that the pattern defined on $[-n;n]^2$
        appears in a configuration of $\sshift{Y}_q$, and is thus $q$-periodic,
        if it appears in a $[-N;N]^2$ pattern not containing any
        $\sshift{Y}_q$-forbidden pattern.

        Without loss of generality, we can assume that the $\sshift{Y}_q-$forbidden
        patterns contain all the patterns defined on $[-n;n]^2$ 
        having distinct values at coordinates separated by $q$ since
        $\sshift{Y}_q$ contains only $q$-periodic configurations.
        We define the SFT $\sshift{Z}$ by forbidding all the patterns defined on
        $[-N;N]^2$ that contain both an $\sshift{X}-$forbidden pattern
        and an $\sshift{Y}_q-$forbidden pattern.
        It is clear that $\sshift{X}\cup\sshift{Y}_q$ is a subset of
        $\sshift{Z}$.

        Let us assume that $\sshift{Z}$ contains a configuration that is not in
        $\sshift{X}\cup\sshift{Y}_q$: $x$ thus contains both an
        $\sshift{X}$ and an $\sshift{Y}_q$-forbidden pattern.
        We show that such an $x$ must contain a pattern defined on $[-N;N]^2$
        containing both an $\sshift{X}$ and an $\sshift{Y}_q-$forbidden pattern
        what contradicts the definition of $\sshift{Z}$.

        Let $i\in\plan$ be where an $\sshift{X}-$forbidden pattern appears in
        $x$. By definition of $\sshift{Z}$, the pattern $x_{|i+[-N;N]^2}$ does
        not contain any $\sshift{Y}_q$-forbidden pattern and is thus 
        $q$-periodic since $\sshift{Y}_q$ forbids all such patterns that are not.
        Therefore, this $\sshift{X}-$forbidden pattern also appears in
        $x$ at positions $i+q\Z$.

        Let $j\in\plan$ be a position where an $\sshift{Y}_q$-forbidden pattern
        appears in $x$.
        The pattern  $x_{|j+[-N;N]^2}$ does not contain, by definition of
        $\sshift{Z}$, any $\sshift{X}-$forbidden pattern; the pattern
        $x_{|j+[-n;n]^2}$ is therefore $p-$periodic for some $p\in\mathcal{P}$
        by the choice of a big enough $N$.
        Since $p$ and $q$ are not collinear, one of $j+p$ or $j-p$ is closer to
        the discrete dotted line $i+q\Z$ than $j$ is and an
        $\sshift{Y}_q$-forbidden pattern appears at these positions too.
        By iterating this reasoning, we obtain a coordinate $j$ where an
        $\sshift{Y}_q$-forbidden pattern appears but for which there exists
        $z\in\Z$
        such that $d(i+qz,j)\leq{}n$: The pattern $x_{|j+[-N;N]^2}$ contains
        both an $\sshift{Y}_q$ and an $\sshift{X}-$forbidden pattern; a
        contradiction to the fact that $x$ belongs to $\sshift{Z}$.
\end{proof}

 With Claim~\ref{claim:perunoncolpersft}, we can make an induction as long as
 $\mathcal{P}$ does not contain two collinear vectors. In order to obtain the
 general case, we group collinear vectors together and thus need to prove that
 this grouping still gives an SFT:
\begin{claim}
        \label{claim:colsft}
    Let $\mathcal{P}$ be a finite subset of $\plan$ such that all the elements
    of $\mathcal{P}$ are collinear and $\sshift{X}$ be an SFT. Then, the
    subshift $\sshift{X}_{\mathcal{P}}\subseteq\sshift{X}$ of configurations of
    $\sshift{X}$ admitting a period in $\mathcal{P}$ is an SFT.
\end{claim}
\begin{proof}[Proof of Claim~\ref{claim:colsft}]
        Let $p=(x,y)$ be a common multiple of all the vectors of
        $\mathcal{P}=\left\{p_1,\ldots,p_n\right\}$.
        Let $\Sigma$ be the alphabet of $\sshift{X}$.
        Consider the one-dimensional subshift $\sshift{Y}$ of
        $(\Sigma^{x\times{}y})^{(x,y)\Z}$ of the configurations that admit a $p_i$ as
        a period when seen as an infinite thick line on $\Sigma^{\Z^2}$.
        Since there is only a finite number of $x\times{}y$-patterns,
        $\sshift{Y}$ is finite and thus an SFT.

        The subshift $\sshift{X}_{\mathcal{P}}$ of $\sshift{X}$ is therefore an
        SFT too: It is defined by the forbidden patterns of $\sshift{X}$ union the
        two-dimensional forbidden patterns of $\sshift{Y}$ in the sense that
        for any $i\in\Z^2$, the configuration defined by
        $x_{|[0;x-1]\times[0;y-1]+(x,y)\Z}$ must belong to $\sshift{Y}$.
\end{proof}

        Let $\mathcal{P}=\left\{p_1,\ldots,p_n\right\}$.
        Let $\mathcal{P}(p_1)$ be the set of vectors of $\mathcal{P}$ that are
        collinear to $p_1$.
        By Claim~\ref{claim:colsft}, $\sshift{X}_{\mathcal{P}(p_1)}$ is an SFT.

        If $\mathcal{P}=\mathcal{P}(p_1)$ then the result is proved.
        Otherwise, let
        $\mathcal{Q}=\mathcal{P}\setminus\mathcal{P}(p_1)=\left\{q_1,\ldots,q_m\right\}$.
        Again, let $\mathcal{Q}(q_1)$ be the set of vectors of $\mathcal{Q}$
        that are collinear to $q_1$ and $q$ be a common multiple of all the
        vectors of $\mathcal{Q}$.
        By Claim~\ref{claim:colsft}, $\sshift{X}_{\mathcal{Q}(q_1)}$ is an SFT.
        By applying Claim~\ref{claim:perunoncolpersft}, with
        $\mathcal{P}=\mathcal{P}(p_1)$,
        $\sshift{X}=\sshift{X}_{\mathcal{P}(p_1)}$, $q=q$ and
        $\sshift{Y}_q=\sshift{X}_{\mathcal{Q}(q_1)}$, we obtain that
        $\sshift{X}_{\mathcal{P}(p_1)}\cup\sshift{X}_{\mathcal{Q}(q_1)}$ is an
        SFT.
        We iterate this reasoning until $\mathcal{P}$ is empty and have proved
        the result.
\end{proof}

\begin{theorem}
        There exists no countable SFT of Cantor-Bendixson rank $\lambda+2$ for
        $\lambda$ a limit ordinal.
\end{theorem}

\begin{proof}
        Suppose there exists such an SFT $\sshift{X}$ of rank $\lambda+2$:
        $\sshift{X}^{(\lambda+2)}=\emptyset$.
        All the configurations of $\sshift{X}^{(\lambda)}$ are at most at
        level~$1$ by Proposition~\ref{prop:pluspetitalorsrkplusgrand} and 
        thus admit a direction of periodicity by Theorem~\ref{thm:sftdnbalorsconf1per}.
        Moreover, by Lemma~\ref{lemma:dersftdifffin}, $\sshift{X}^{(\lambda)}$
        is finite up to translations:
        Let $\mathcal{P}\subseteq{}\Z^2$ be a finite set of periods such that any
        configuration of $\sshift{X}^{(\lambda)}$ admits a period from
        $\mathcal{P}$.
        By applying Lemma~\ref{lemma:perfinsft}, we know that the set
        configurations of $\sshift{X}$ admitting a period from $\mathcal{P}$ is
        an SFT.
        Let $\sshift{X}_{\mathcal{P}}$ denote this SFT and $\ForbPat$ the
        forbidden patterns defining it.
        Let $\sshift{X}_{\left\{p\right\}}$ denote the configurations of
        $\sshift{X}$ admitting $p$ as a period.
        
        If for a $\beta < \lambda$, all the configurations of
        $\sshift{X}^{(\beta)}$ contain only configurations of
        $\sshift{X}_{\mathcal{P}}$,
        then $\sshift{X}^{(\beta)}$ is included in 
        $\cup_{p\in \mathcal{P}} \sshift{X}_{\left\{p\right\}}$.
        $\sshift{X}_{\left\{p\right\}}$ is a countable SFT that can be seen as
        a one-dimensional countable SFT (and then repeated periodically along
        the $p$ direction), which thus has finite Cantor-Bendixson rank by \eg
        \cite[Proposition~$3.8$]{villejournalcb}.
        This would mean that $\sshift{X} \cap \sshift{X}_{\left\{p\right\}}$
        would have finite Cantor-Bendixson rank, and thus $\sshift{X}^{(\beta)}$
        too as it would be included in the finite union of closed sets with
        finite Cantor-Bendixson rank.

        Since we assumed $\lambda$ infinite, for any $\beta < \lambda$,
        $\sshift{X}^{(\beta)}$ contains a configuration that is not in
        $\sshift{X}_{\mathcal{P}}$ and thus which contains a pattern in
        $\ForbPat$.
        Now, $\lambda = \cup \beta_i$ ; take a sequence of $x_i$'s isolated in
        $\sshift{X}^{(\beta_i)}$ that contains a pattern of $\ForbPat$ at its
        center.
        By compactness, one can extract a converging subsequence of these
        $x_i$'s whose limit is in $\sshift{X}^{(\lambda)}$ but this limit would
        contain a pattern of $\ForbPat$ at its center: A contradiction with the
        fact that we assumed $\sshift{X}^{(\lambda)}$ to contain only
        configurations with a period from $\mathcal{P}$.
\end{proof}

\section*{Conclusions and future work}
\label{sec:concl}
As a short summary of the work presented in this paper, we applied two different
notions to the study of subshifts: the Cantor-Bendixson rank which comes from
topology and the pre-order $\preceq$ which is more combinatorial.
These two notions allowed us to study in great details the structure of
subshifts: We characterized the configurations in finite
subshifts~\cite{stacsstructuralaspect}, we have been
able to capture exactly the causes of the uncountability of subshifts
(Theorem~\ref{thm:caracsshiftnondnb}) and, finally, the major part of this paper
is devoted to describing the structure of countable subshifts.

We proved the best possible existence results on the levels for $\preceq$:
Level~$0$ and $1$ always exist and level~$2$ exists in countable SFTs. Those
results are optimal since level~$2$ may not exist in countable sofic shifts and
so does level~$3$ in countable SFTs~\cite{villejournalcb}.
While configurations at level~$0$ are easy to characterize, we characterized
those at level~$1$ in countable subshifts. The only missing characterization is
for configurations at level~$2$ in countable SFTs.

An interesting connection is with the theory of the complexity
function. The complexity function $\complf{m}{n}{c}$ of a
configuration $\conf{c}$ is the number of distinct  $m\times{}n$ patterns that appear in $\conf{c}$.
It is obvious that if $\conf{c} \preceq \conf{c'}$ then, for all $m$ and $n$,
$\complf{m}{n}{c} \leq \complf{m}{n}{c'}$.
The connection seems however deeper: In all our SFT examples,
the minimal elements have a \emph{constant} complexity function,
the elements at level $1$ a \emph{linear} complexity function,
and the elements at level $2$ a \emph{quadratic} complexity function.
We do not know if this statement is a coincidence or part of a bigger theorem.
In particular, if $\conf{c} \prec \conf{c'}$ then is 
$\complf{m}{n}{c}$ a $o(\complf{m}{n}{c'})$ ?
This connection, if it exists, might help to investigate the 
Nivat conjecture~\cite{nivatconj}: If
$\complf{m}{n}{c}\leq{}mn$ for some $m,n$ and some two-dimensional configuration
$\conf{c}$, then $\conf{c}$ admits one direction of periodicity.
For example, the conjecture is true if $c$ is such that $\engendre{c}$ is a
countable SFT of Cantor-Bendixson rank at most $3$ by Theorem~\ref{thm:rang2}.

\bibliographystyle{plain}
\bibliography{biblio/article,biblio/ca-faq,biblio/stacs08,biblio/latin10}
\end{document}